\documentclass[a4,12pt]{article}
\usepackage{fullpage}
\usepackage{amsfonts}
\usepackage{amsmath}
\usepackage{amssymb} 
\usepackage{graphicx}
\usepackage{xcolor}
\usepackage{hyperref}
\usepackage[mathscr]{eucal}
\newcommand{\R}{\mathbb R}

\title{Slow-fast systems with an equilibrium near the folded slow manifold}
\author{N. G. Gelfreikh, A. V. Ivanov}
\date{}
\begin{document}
\renewcommand{\theequation}{\arabic{section}.\arabic{equation}}
\maketitle

\begin{abstract}
We study a slow-fast system with two slow and one fast variables. 
We  assume that the slow manifold of the system possesses a fold and there is  an equilibrium of the system 
in a small neighbourhood of the fold. We derive a normal form for the system 
in a neighbourhood of the pair "equilibrium-fold" 
and study the dynamics of the normal form. In particular, as the ratio of two time scales tends to zero we obtain an asymptotic formula for the Poincar\'e map
and calculate the parameter values for the first period-doubling bifurcation. The theory is applied to a generalization
of the  FitzHugh-Nagumo system.
\end{abstract}



\section{Introduction}

Since fundamental works of A. N. Tikhonov, L. S. Pontryagin, N. Fenichel \cite{Tikh}, \cite{Pon}, \cite{Fen} singular perturbed dynamical systems became a subject of many researches and are intensively studied nowadays. These systems have many applications in various areas of physics: mechanics, hydrodynamics, plasma physics, neurobiology and others. They are aimed to describe those systems for which different processes are running with different speeds. Mathematically, such systems can often be written in the form of so-called slow-fast system
\begin{equation}\label{eqn_slowfast_1}
\begin{aligned}
\varepsilon \dot x &=F(x,y) ,\\
\dot y &= G(x,y),
\end{aligned}
\end{equation}
where $x\in\R^n$ is a fast variable and $y\in\R^m$ is a slow one. The parameter $\varepsilon$ describes the ratio of two time scales and is assumed to be small. 
In the limit $\varepsilon = 0$ we obtain the slow system
\begin{equation}\label{eqn_slow}
\begin{aligned}
0&=F(x,y) ,\\
\dot y &= G(x,y),
\end{aligned}
\end{equation}
which describes the motion in a vicinity of the slow manifold 
defined by the equality $F(x, y) = 0$. On the other hand, using the fast time $s = t/\varepsilon$ 
we can rewrite the system (\ref{eqn_slowfast_1}) as
\begin{equation}\label{eqn_slowfast_2}
\begin{aligned}
x' &=&F(x,y) ,\\
y' &=& \varepsilon G(x,y),
\end{aligned}
\end{equation}
where the prime stands for the derivative with respect to $s$.
Setting $\varepsilon = 0$ in (\ref{eqn_slowfast_2}), we obtain the fast system
\begin{equation}\label{eqn_fast}
\begin{aligned}
x'&=F(x,y),\\
y'&= 0.
\end{aligned}
\end{equation}
 The fast system approximates the original system on any finite interval with respect to $s$
due to smooth dependence of solutions on a vectorfield. It should be noted that 
for $\varepsilon\neq 0$ the systems (\ref{eqn_slowfast_1})
and (\ref{eqn_slowfast_2}) are equivalent, 
however, the limit systems (\ref{eqn_slow}) and (\ref{eqn_fast}) are different. 

The standard tool for studying  slow-fast systems 
is based on the geometric singular perturbation theory (GSPT) \cite{Fen}.
Using the notion of normal hyperbolicity,
GSPT predicts that a trajectory attracted by a stable branch of the slow manifold follows closely a trajectory of (\ref{eqn_slow})
till this trajectory hits a singularity of the slow manifold. 
Singularities of slow manifolds cause a variety of phenomena
including  delay of stability loss and canard explosion
(see e.g. \cite{Nei87}, \cite{Nei88}, \cite{Szm}, \cite{KruSzm}).

The present paper was inspired by \cite{Zaks}, where the author studied a  FitzHugh-Nagumo-like sysytem originated from the mathematical theory of neural cells. This system consists of three ODEs with one fast variable corresponding to the membrane potential and two slow gating variables:
\begin{eqnarray}
\nonumber
\varepsilon \dot x&=&x-x^3/3-y-z ,\\
\label{eq_FHN}
\dot y&=& a+x ,\\
\nonumber
\dot z &=& a+x-z,
\end{eqnarray}
where $\varepsilon$ is a small parameter and $a$ is a real parameter.
 The slow manifold of the system is described by the equation $x-x^3/3-y-z = 0$ and possesses folds at $x = \pm 1, y+z = \pm 2/3$.
The system has a unique equilibrium which is
close to the fold if $a$ is close to one.
The equilibrium  is stable for larger values of $a$
and undergoes a supercritical Andronov-Hopf bifurcation at  $a_{H} = 1-\frac{1}{4}\varepsilon + O(\varepsilon^{2})$ (see \cite{Zaks} for details) .

In \cite{Zaks} M. Zaks found numerically that the initial periodic orbit may lose stability via a sequence of period-doubling bifurcations.  Studying numerically the period-doubling cascades for small but fixed values of the parameter $\varepsilon$, M. Zaks observed that the cascade follows the Feigenbaum law with the Feigenbaum constant $4.67\ldots$ typical for dissipative systems. On the other hand
for smaller values of $\varepsilon$ the process switches to the Feigenbaum constant of a conservative map as, in the limit $\varepsilon\to 0$,
two-dimensional  Poincar\'e map nearly preserves the area. The reason for such phenomenon was assumed to be in the closeness of the equilibrium to a fold of the slow manifold. 

In the present paper we consider a family of slow-fast systems having one fast and two slow variables and depending on a real parameter $\delta $:
\begin{eqnarray}
\nonumber
\varepsilon \dot x&=& F(x,y,z;\delta )  \,,\\
\label{Eq:initialODE_0}
\dot y&=&  G_1(x,y,z;\delta )  \,,\\
\nonumber
\dot z &=&  G_2(x,y,z;\delta )  \, ,
\end{eqnarray}
where $(x,y,z)\in \mathbb{R}^{3}$,  functions $F$, $G_1$ and $G_2$ are smooth functions, $\varepsilon$ and $\delta $ are small parameters. 

 We suppose that the system (\ref{Eq:initialODE_0}) possesses an equilibrium $(x,y,z) = (x_0(\delta ), y_0(\delta ), z_0(\delta ) )$
and if $\delta =0$ the equilibrium lies on the fold of the slow manifold.
Shifting the origin into the equilibrium, we have the following conditions:
\begin{equation}
  \label{eq_Cond_RHS1}
  F(0,0,0;\delta )=0, \  G_1(0,0,0;\delta )=0, \ G_2(0,0,0;\delta )=0.
\end{equation}  

The slow manifold is defined by the equality
\begin{equation*}
F(x,y,z;\delta ) = 0.
\end{equation*}
We assume that for $ \delta =0$ the slow manifold possesses a non-degenerate fold which is tangent to a fast fibre. 
More precisely, 
we impose the following conditions (see e.g. \cite{Kuehn}):
\begin{equation}
\label{eq_Cond_RHS2}
 F'_x(0,0,0;0)=0,\  \nabla F(0,0,0;0)\neq 0, \  F''_{x^2}(0,0,0;0)\ne 0.
\end{equation}

Finally, we impose a condition
\begin{equation}
\label{eq_Cond_RHS}
 F'_y (0,0,0; 0) G'_{1x}(0,0,0; 0) +F'_z (0,0,0; 0) G'_{2x} (0,0,0; 0) < 0,
\end{equation}
which ensures the linear stability of the equilibrium in the following sence. Consider the equations (\ref{Eq:initialODE_0}) linearized at the point $(0, 0, 0, 0)$. Then its characteristic equation reads
\begin{multline*}
\label{Eq:charact}
-\lambda^{3}+(G'_{1y}+G'_{2z})\lambda^{2} + \left(\varepsilon^{-1}(F'_{y}G'_{1x} + F'_{z}G'_{2x}) - (G'_{1y}G'_{2z}-G'_{1z}G'_{2y})\right)\lambda - \\
- \varepsilon^{-1}\left(F'_{y}(G'_{1x}G'_{2z}-G'_{1z}G'_{2x}) - F'_{z}(G'_{1x}G'_{2y}-G'_{1y}G'_{2x})\right) = 0,
\end{multline*}
where all the derivatives are evaluated at the point $(0,0,0,0)$. Substituting $\Lambda = \sqrt{\varepsilon}\lambda$, we obtain
\begin{equation}
\label{Eq:charact_mod}
-\Lambda\left(\Lambda^{2} - (F'_{y}G'_{1x} + F'_{z}G'_{2x})\right) = O(\sqrt{\varepsilon}).
\end{equation}
By setting $\varepsilon=0$ in (\ref{Eq:charact_mod}) one arrives at the limit equation
\begin{equation*}
-\Lambda\left(\Lambda^{2} - (F'_{y}G'_{1x} + F'_{z}G'_{2x})\right) = 0,
\end{equation*}
which has three roots $\Lambda_{0} = 0$, $\Lambda_{\pm} = \pm\sqrt{F'_{y}G'_{1x} + F'_{z}G'_{2x}}$. Then the condition (\ref{eq_Cond_RHS}) guarantees that all these roots have non-positive (in fact, zero) real parts.

Under these assumptions we derive a  normal form 
for  the system 
\eqref{Eq:initialODE_0}
in a neighbourhood of the pair "fold-equilibrium":
\begin{eqnarray*}
\xi'&=&\xi^2-\eta +
 \mu (\gamma_0 \xi +\gamma \xi^3) + 
 \mu^2 g_1(\xi, \eta, \zeta)  \,, \\
\eta'&=&  2\xi + \mu (\alpha_1 \eta +  \alpha_2 \zeta )
 +    \mu^2 g_2(\xi, \eta, \zeta)  \,, \\
 \zeta' &=&   \mu ( \beta_1 \eta +  \beta_2 \zeta )
 +   \mu^2 g_3(\xi, \eta, \zeta),  \, 
\end{eqnarray*}
where $\mu = \sqrt{\varepsilon}$ is a new small parameter, $\gamma_{0}, \gamma, \alpha_{1,2}, \beta_{1,2}$ are constants and functions $g_{k}, k=1,2,3$ are polynomials of $(\xi, \eta, \zeta)$, which will be specified later.  The parameter $\gamma_{0}$ has a special role: it describes the closeness of the equilibrium to the fold. Using the Poincar\'e map technique we study the dynamics of the normal form. In particular, we show that in a $\mu\ln \mu$-neighborhood of the equilibrium the normal form system has a periodic trajectory. Varying distance between the equilibrium and the fold one may observe how this thrajectory undergoes the period-doubling bifurcation. We obtain conditions on parameters of the normal form which correspond to this scenario.
 
The paper is organized in the following way.
In the second section we derive the normal form.
The asymptotics for the Poincar\'e map 
associated with  the normal form are obtained in Section 3.
Section 4 is devoted to a construction of asymptotic
conditions for existence of a periodic orbit and its period-doubling bifurcations. The main result of the paper is the asymptotic formula (\ref{eq_FPD}). Finally, in Section 5 we apply the obtained results to the FitzHugh-Nagumo system and 
compare them with numeric data.

\section{Normal form}

In this section we derive the formal normal form for the system (\ref{Eq:initialODE_0}) when the right-hand side satisfies the assumptions (\ref{eq_Cond_RHS1}), (\ref{eq_Cond_RHS2}) and (\ref{eq_Cond_RHS}).

We construct  a sequence of changes of space variables and time and apply them  to (\ref{Eq:initialODE_0}).

%
%
%
%
%
%
First, we introduce a new small parameter $\mu$:
\[ \varepsilon =\mu^2 
\]
and make a rescaling of the space variables, time and the second parameter $\delta$:
\begin{equation}\label{eq_CoV_1}
x=\mu X_{1},\quad  y=  \mu^2 Y_{1},\quad
z =  \mu^2 Z_{1},\quad t =  \mu s,\quad
\delta = \mu^{2}\sigma  \, . 
\end{equation}
Then 
\(
\varepsilon \dot x=\mu^2 X'_{1}   \,,\
\dot y=\mu Y'_{1} \, ,\
\dot z = \mu Z'_{1}   \, ,
\)
where the prime stands for the derivative with respect to the "semi-fast" time $s$.
We substitute (\ref{eq_CoV_1}) into (\ref{Eq:initialODE_0}) and use the Taylor formula for the right hand side of (\ref{Eq:initialODE_0}) in a neighborhood of the point $(x,y,z;\delta ) = (0,0,0;0)$. Then, taking into account (\ref{eq_Cond_RHS1}) and (\ref{eq_Cond_RHS2}), one obtains
\begin{eqnarray}
\label{Eq:initialODE_3}\nonumber
 X'_{1}&=&  F'_y Y_{1} +F'_z Z_{1} + \frac12 F''_{x^2} X_{1}^2 + 
 \\ & & +\mu \left[ \sigma F''_{x\delta } X_{1} +
 F''_{xy} X_{1}Y_{1} + F''_{xz} X_{1}Z_{1} +\frac16 F'''_{x^3} X_{1}^3
 \right] + \\
 \nonumber & &+  \mu^2 \left[\sigma F''_{y\delta }Y_{1} + \sigma F''_{z\delta } Z_{1} + \frac12 F''_{y^2}Y_{1}^2 +\frac12 F''_{z^2} Z_{1}^2 +F''_{yz} Y_{1}Z_{1} 
\right. + \\ \nonumber  & &+ \left. 
 X_{1}^2 \left( \frac12 \sigma F'''_{x^2\delta } +
 \frac12 F'''_{x^2y}Y_{1} + \frac12 F'''_{x^2z}Z_{1}  \right) + \frac{1}{24} F^{(4)}_{x^4} X_{1}^4
  \right] +O(\mu^3),
 \\
\nonumber
 Y'_{1}&=&  G'_{1x} X_{1} + \mu \left[  G'_{1y}Y_{1} +  G'_{1z}Z_{1} + \frac12  G''_{1x^2} X_{1}^2 \right] + \\  \nonumber & &
+\mu^2
 \left[ X_{1} \left(  \sigma G''_{1x\delta } +
  G''_{1xy}Y_{1} +  G''_{1xz} Z_{1} \right) + \frac16  G'''_{1x^3}X_{1}^3 \right] +O(\mu^3)  \,, \\
\nonumber
  Z'_{1}&=& G'_{2x} X_{1} + \mu \left[  G'_{2y}Y_{1} + G'_{2z}Z_{1} + \frac12  G''_{2x^2} X_{1}^2 \right] + \\ \nonumber & &
+\mu^2
 \left[ X_{1} \left(  \sigma G''_{2x\delta } +
   G''_{2xy}Y_{1} +  G''_{2xz} Z_{1} \right) + \frac16  G'''_{2x^3}X_{1}^3 \right] +O(\mu^3)   \, .
\end{eqnarray}
Here and below in this section all derivatives of the functions $F$, $G_1$ and $G_2$ are evaluated at the point $(x,y,z; \delta )=(0,0,0;0)$. 

For $\mu=0$ the system takes the form
\begin{eqnarray*}
X'_{1}&=&  F'_y Y_{1} +F'_z Z_{1} + \frac12 F''_{x^2} X_{1}^2 ,\\
Y'_{1}&=&  G'_{1x} X_{1} ,\\
  Z'_{1}&=& G'_{2x} X_{1} .
\end{eqnarray*}
 and the multipliers of the corresponding linearized system satisfy an equation 
\[ - \lambda^3 + \lambda (F'_y G_{1x}' + F'_z G_{2x}') =0. \]
As it was mentioned above, in this paper we consider the case of stable equilibrium. Then, denoting by
\begin{equation}\label{def_D}
D=F'_y G_{1x}' + F'_z G_{2x}',
\end{equation} 
the condition  (\ref{eq_Cond_RHS}) takes the form
\begin{equation}\label{eq_D} 
D<0. 
\end{equation}

%
%
%
%
%
One may note that for $\mu=0$ the system (\ref{Eq:initialODE_3}) has an obvious integral:
\[ - G_{2x}' Y_{1} + G_{1x}' Z_{1}.\]
Taking this into account, we introduce new variables $(X_{2}, Y_{2}, Z_{2})$ by
\[X_{2} = X_{1}, \qquad Y_{2} = F'_y Y_{1}+F'_z Z_{1}, 
\qquad Z_{2} = - G_{2x}' Y_{1} + G_{1x}' Z_{1}.
\]
Due to (\ref{eq_D}) the inverse change of variables is well-defined:
\[ X_{1} = X_{2}, \qquad Y_{1} = \frac{1}{D} ( G_{1x}'  Y_{2} - F'_z Z_{2} ), \qquad Z_{1} =  \frac{1}{D} (G_{2x}' Y_{2} + F'_y Z_{2}) \]
and the system (\ref{Eq:initialODE_3}) can be rewritten as
\begin{eqnarray}
\label{Eq:initialODE_4}
X'_{2}&=&  Y_{2} + \frac12 F''_{x^2} X_{2}^2 +\mu \left\{  X_{2}
    \left[  \gamma_0^{(2)} +  \gamma_1^{(2)} Y_{2}
    +  \gamma_2^{(2)} Z_{2} \right] +  \gamma_3^{(2)} X_{2}^3 \right\}+
 \\ \nonumber &  & \mu^2 \biggl\{    \gamma_4^{(2)} Y_{2}
 + \gamma_5^{(2)} Z_{2} +
 \gamma_6^{(2)} Y_{2}^2 
 +  \gamma_7^{(2)} Z_{2}^2   
 +  \gamma_{8}^{(2)} Y_{2} Z_{2} +
\\ \nonumber &  & X_{2}^2 \left[  \gamma_{9}^{(2)} + \gamma_{10}^{(2)} Y_{2}  
 +   \gamma_{11}^{(2)} Z_{2} 
 \right] +   \gamma_{12}^{(2)} X_{2}^4 \biggr \}
  +O(\mu^3),\\
 \nonumber
Y'_{2}&=& D X_{2} +
\mu \left(  \alpha_1^{(2)} Y_{2} 
+  \alpha_2^{(2)} Z_{2}
+  \alpha_3^{(2)} X_{2}^2 \right) +
 \mu^2  \left\{
X_{2} \left[    \alpha_4^{(2)} + \alpha_5^{(2)} Y_{2} +  \alpha_6^{(2)}  Z_{2}  \right]
+ \alpha_7^{(2)} X_{2}^3  \right\}
 +O(\mu^3)  \,, \\
 \nonumber
Z'_{2}&=& \mu \left[  \beta_1^{(2)} Y_{2} +
 \beta_2^{(2)}
 Z_{2} +  \beta_3^{(2)} X_{2}^2  \right] + \mu^2 \left\{  X_{2}  \left[   \beta_4^{(2)} + \beta_5^{(2)} Y_{2} +  
\beta_6^{(2)} Z_{2} 
 \right] +  \beta_7^{(2)} X_{2}^3   \right\} +O(\mu^3)  \, ,
\end{eqnarray}
where
\begin{equation} \label{eq_Coef2}
  \gamma_0^{(2)} = \sigma F''_{x\delta } ,\quad
 \gamma_1^{(2)} = \frac{1}{D} \left( F''_{xy} G'_{1x}
    + F''_{xz} G'_{2x} \right) ,\quad  \gamma_2^{(2)} =
\frac{1}{D} \left(F''_{xz} F'_y - F''_{xy} F'_z \right), \quad
\gamma_3^{(2)} = \frac16 F'''_{x^3} ,
\end{equation}

\[  \alpha_1^{(2)} = \frac{1}{D} (F'_yG'_{1x}G'_{1y} + F'_y G'_{1z} G'_{2x} +  F'_zG'_{1x} G'_{2y} +  F'_z G'_{2x} G'_{2z} ), \]
\[  \alpha_2^{(2)} = \frac{1}{D} (F'^2_y G'_{1z} - F'^2_z G'_{2y}
+F'_yF'_z G'_{2z} - F'_yF'_z G'_{1y} ),
\quad   \alpha_3^{(2)} = \frac12 (F'_y G''_{1x^2} + F'_z G''_{2x^2}), \]

\[  \beta_1^{(2)} = \frac{1}{D} (G'^2_{1x} G'_{2y}
-  G'_{1x} G'_{2x} G'_{1y} + G'_{1x} G'_{2x} G'_{2z}
- G'^2_{2x} G'_{1z}),\]
\[  \beta_2^{(2)} = \frac{1}{D} ( - F'_zG'_{1x}G'_{2y}
+ F'_zG'_{2x}G'_{1y} + F'_y G'_{1x}G'_{2z} - F'_y G'_{2x} G'_{1z}), \]
\[ \beta_3^{(2)} = \frac12 (G'_{1x}G''_{2x^2} - G'_{2x} G''_{1x^2}).\]
We do not write formulae for the coefficients of terms of the order $\mu^2$ (i.e. $\gamma^{(2)}_i$, $\alpha^{(2)}_i$, $\beta^{(2)}_i$ with $i \ge 4$) since they do not enter the main result of the paper.
\newline

In the system (\ref{Eq:initialODE_4}) the variable $Z_{2}$ is slow and the leading term of the right-hand side in the first equation does not contain $Z_{2}$.

%
%
%
%
%
%

The next change of variables is aimed at simplification of the leading term of the right hand side in (\ref{Eq:initialODE_4}). Using scaling
\[ X_{2}=k X_{3}, \quad Y_{2}=m Y_{3},\quad Z_{2}=Z_{3},\quad s =n \tau\]
we represent (\ref{Eq:initialODE_4}) as
\begin{eqnarray}
\nonumber
X'_{3}&=&\frac{nm}{k}Y_{3}+\frac12 F''_{ x^2} kn X_{3}^2+O(\mu), \\
\nonumber
Y'_{3}&=& D \frac{kn}{m} X_{3}+O(\mu)  \,, \\
\nonumber
Z'_{3} &=& O(\mu)  \, .
\end{eqnarray}
We fix the factors $m, k, n$ in the following way
\[ m = \frac{D}{F''_{x^2}}, \quad k =  \frac{\sqrt{-2D}}{F''_{x^2}}  , \quad n = \sqrt{ \frac{2}{-D}}\, .  
\]
Then the leading term of (\ref{Eq:initialODE_4}) is simplified to 
\begin{eqnarray}
\nonumber
X'_{3}&=&X_{3}^2-Y_{3}\,, \\
\nonumber
Y'_{3}&=& 2X_{3}\,, \\
\nonumber
Z'_{3} &=& 0   \, . 
\end{eqnarray}
Here the prime stands for the derivative with respect to $\tau $.
And the whole system (\ref{Eq:initialODE_4}) takes the form
\begin{eqnarray}
\label{Eq:initialODE_5}
X'_{3}&=&X_{3}^2-Y_{3}
+\mu \left( X_{3}(\gamma_0^{(3)} + \gamma_1^{(3)}Y_{3}+\gamma_2^{(3)}Z_{3})+\gamma_3^{(3)}X_{3}^3\right) +
\\ \nonumber
&+& \mu^2 \biggl\{    \gamma_4^{(3)} Y_{3}
 + \gamma_5^{(3)} Z_{3} +
 \gamma_6^{(3)} Y_{3}^2 
 +  \gamma_7^{(3)} Z_{3}^2   
 +  \gamma_{8}^{(3)} Y_{3} Z_{3} + 
\\ \nonumber
&+& X_{3}^2 \left[  \gamma_{9}^{(3)} + \gamma_{10}^{(3)} Y_{3}  
 +   \gamma_{11}^{(3)} Z_{3} 
 \right] +   \gamma_{12}^{(3)}X_{3}^4 \biggr \}+
O(\mu^3)\,, \\
\nonumber
Y'_{3}&=&  2X_{3}+\mu (\alpha_1^{(3)}Y_{3}+ \alpha_2^{(3)} Z_{3}+  \alpha_3^{(3)} X_{3}^2 ) + \\ \nonumber
&+& \mu^2 ( X_{3}(  \alpha_4^{(3)} + \alpha_5^{(3)} Y_{3}+ \alpha_6^{(3)}Z_{3}) +  \alpha_7^{(3)}X_{3}^3 )+O(\mu^3)  \,, \\
\nonumber
Z'_{3} &=&  \mu \left(  \beta_1^{(3)} Y_{3}+  \beta_2^{(3)} Z_{3} +  \beta_3^{(3)} X_{3}^2 \right)
 +   \\ \nonumber
&+& \mu^2 (X_{3}(  \beta_4^{(3)} + \beta_5^{(3)} Y_{3}+ \beta_6^{(3)} Z_{3}) + \beta_7^{(3)} X_{3}^3 )+O(\mu^3)   \, ,
\end{eqnarray}
where
\begin{equation} \label{eq_Coef3}
\gamma_0^{(3)} = \sqrt{ \frac{2}{-D}} \,  \gamma_0^{(2)} ,\quad  
\gamma_1^{(3)}= -\frac{\sqrt{-2D}}{F''_{x^2}} \,  \gamma_1^{(2)},\quad 
\gamma_2^{(3)}= \sqrt{ \frac{2}{-D}} \, \gamma_2^{(2)},\quad 
\gamma_3^{(3)}= \frac{2\sqrt{-2D}}{F''^2_{x^2}}\,  \gamma_3^{(2)} , 
\end{equation}

\[ \alpha_1^{(3)} = \sqrt{ \frac{2}{-D}} \,  \alpha_1^{(2)} ,\quad  
\alpha_2^{(3)} = \sqrt{ \frac{2}{-D}} \,    
\frac{ F''_{x^2}}{D}\,  \alpha_2^{(2)} ,\quad 
\alpha_3^{(3)} = -\frac{2\sqrt{2}}{F''_{x^2}\sqrt{-D}}\,  \alpha_3^{(2)} ,
\]

\[ \beta_1^{(3)} = -\frac{\sqrt{-2D}}{F''_{x^2}}\,  \beta_1^{(2)}, \quad
\beta_2^{(3)} = \sqrt{ \frac{2}{-D}} \,   \beta_2^{(2)}, \quad
\beta_3^{(3)} = \frac{2\sqrt{-2D}}{F''^2_{x^2}}\,  \beta_3^{(2)}. \]
%
%
%
%
%
%
The final change of variables is aimed at excluding as many coefficients $\alpha_{i}, \beta_{i}, \gamma_{i}$ as possible. One may remark that equations (\ref{Eq:initialODE_5}) (and in fact already (\ref{Eq:initialODE_3}) due to scaling (\ref{eq_CoV_1})) possess a symmetry. Namely, they are invariant with respect to the following transformation
\begin{equation}\label{eq_sym}
(X_{3},\tau,\mu) \mapsto (-X_{3}, -\tau, -\mu)
\end{equation}
This symmetry will be used to simplify the study of the Poincar\'e map.  Thus, our purpose is not only to exclude a number of coefficients and keep the leading term, but also to preserve this symmetry.

For this reason we consider the following close-to-identity change of variables 
\begin{equation}\label{eq_CoV_4}
\left(
\begin{array}{c}
X_{3}\\
Y_{3}\\
Z_{3}
\end{array}\right)=(I+\mu A+\mu^2 B)\left(
\begin{array}{c}
\xi \\
\eta \\
\zeta
\end{array}\right) ,
\end{equation}
where
\[ A=\left(\begin{array}{ccc}
0 & B_1 & C_1 \\
A_2 & 0 & 0\\
A_3 & 0 & 0
\end{array}
\right)
,\qquad B=\left(\begin{array}{ccc}
a_1 & 0 & 0 \\
0 & b_2 & c_2\\
0 & b_3 & c_3
\end{array}
\right)
. \]
Inverting of (\ref{eq_CoV_4}) yields
\begin{equation}\label{eq_CoV_4_inv}
\left(\begin{array}{c}
\xi\\
\eta\\
\zeta
\end{array}
\right)=(I-\mu A+\mu^2 (A^2- B))\left(\begin{array}{c}
X_{3}\\
Y_{3}\\
Z_{3}
\end{array}
\right)+O(\mu^3) .
\end{equation}

In order to simplify the terms of the first order, we choose
\[ B_1 = -\frac12 \gamma_1^{(3)},\ C_1 = - \frac12 \gamma_2^{(3)}, \ A_2 = \alpha_3^{(3)} ,\  A_3 = \beta_3^{(3)} .\]
Then by appropriate choice of $a_1$, $b_2$, $b_3$ and $c_2$ we remove four coefficients of the order $\mu^2$ and
 obtain the following system:
\begin{eqnarray*}
\xi'&=&\xi^2-\eta +
 \mu f_1+ 
 \mu^2 g_1 +O( \mu^3) \,, \\
\eta'&=&  2\xi + \mu f_2
 +    \mu^2 g_2 +O( \mu^3) \,, \\
 \zeta' &=&   \mu f_3
 +   \mu^2 g_3 +O( \mu^3)  \, ,
\end{eqnarray*}
where 
\begin{eqnarray}
\nonumber
  f_1 &=& \gamma_0 \xi +\gamma \xi^3 ,\\
\nonumber
  g_1 &=&   \gamma_1 \eta  + \gamma_2  \eta^2+ \gamma_3 \zeta^2
+   \gamma_4\eta\zeta + \xi^2  ( \gamma_5\eta+  \gamma_6\zeta ) +  \gamma_7\xi^4, \\
\nonumber
f_2 &=&  \alpha_1 \eta +  \alpha_2 \zeta  ,\\
\label{eq_rhs_nf}
g_2 &=&  \xi ( \alpha_3\eta + \alpha_4\zeta  ) +  \alpha_5  \xi^3 ,\\
\nonumber
f_3 &=&  \beta_1 \eta +  \beta_2 \zeta , \\
\nonumber
g_3 &=& \xi ( \beta_3 \eta + \beta_4 \zeta ) +\beta_5 \xi^3 
\end{eqnarray}
and
\begin{equation} \label{eq_Coef}
\gamma_0 = \gamma_0^{(3)} +\gamma_1^{(3)}-\alpha_3^{(3)} ,\qquad \gamma = \gamma_3^{(3)},
\end{equation}
\[ \alpha_1 = \alpha_1^{(3)} + \alpha_3^{(3)} - \gamma_1^{(3)}, \quad
\alpha_2 = \alpha_2^{(3)} -\gamma_2^{(3)}, \quad
\beta_1 = \beta_1^{(3)} + \beta_3^{(3)}, \quad
\beta_2 = \beta_2^{(3)}. \]

We omit terms of the order $O(\mu^{3})$ and obtain the following system:
\begin{eqnarray}
\nonumber \label{eq_NormForm}
\xi'&=&\xi^2-\eta +
 \mu f_1+ 
 \mu^2 g_1  \,, \\
\eta'&=&  2\xi + \mu f_2
 +    \mu^2 g_2  \,, \\
\nonumber
 \zeta' &=&   \mu f_3
 +   \mu^2 g_3   \, 
\end{eqnarray}
with $f_i$ and $g_i$ defined by (\ref{eq_rhs_nf}). 
The system (\ref{eq_NormForm}) will be called the normal form in a neighborhood of a pair "equilibrium-fold". One should emphasize that equations (\ref{eq_NormForm}) are invariant under the following transformation
\begin{equation}\label{eq_Symm}
(\xi, \tau, \mu) \mapsto (-\xi, -\tau, -\mu).
\end{equation}

\section{Dynamics of the normal form}

\subsection{Poincar\'e map}

\begin{figure}[h] 
\begin{center}
\includegraphics[width=8cm]{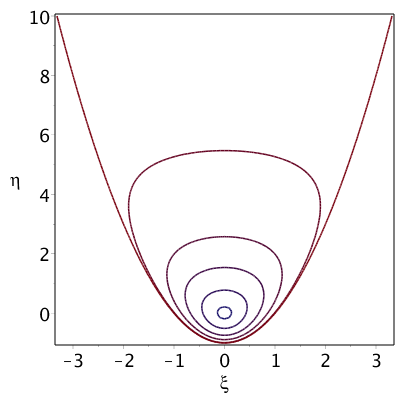}
\end{center}
\caption{Orbits of the unperturbed system (\ref{eq_Unpert_v}) for different values of $J_0$.} 
\label{Fig:unpert_orbits}
\end{figure}

Setting $\mu =0$ in (\ref{eq_NormForm}) (that corresponds to $\varepsilon =0$ for the original system), we obtain the  system
\begin{eqnarray}
\nonumber
\xi'&=&\xi^2-\eta, \\
\label{eq_Unpert_v}
\eta'&=& 2\xi  \,, \\
\nonumber
\zeta' &=& 0  \, .
\end{eqnarray}
In addition to the obvious integral of motion $\zeta$, the unperturbed system has the second integral
\begin{equation}\label{eq_J}
J=(\eta + 1 - \xi^2)e^{-(\eta+1)}.
\end{equation}
The orbits of the system (\ref{eq_Unpert_v}) belong to  
intersections of the integrals' level sets (see fig.\ref{Fig:unpert_orbits})
\begin{equation}\label{eq_unpertOrbit}
\zeta = \zeta_0,\quad (\eta + 1 - \xi^2)e^{-(\eta+1)} = J_0.
\end{equation}
Note that all fixed points of (\ref{eq_Unpert_v}) belong to 
the line $(0, 0, \zeta)$. Moreover, the unperturbed system possesses a separatrix. 
Namely, the parabola $\eta = \xi^2-1$ (that corresponds to $J_0=0$) separates the plane $(\xi,\eta)$ into two parts:
above the parabola ($J_0>0$) all orbits of (\ref{eq_Unpert_v}) are closed and below it ($J_0<0$) all orbits are not closed.

To study the dynamics of (\ref{eq_NormForm}) we construct the Poincar\'e map for this system. The normal form system can be considered as a small perturbation of (\ref{eq_Unpert_v}). 

We add to the system (\ref{eq_NormForm}) an additional equation in order to describe evolution of the variable~$J$. Using (\ref{eq_J}) and (\ref{eq_NormForm}), we obtain the following equation on $J$:
\begin{equation}\label{eq_dynJ}
J' = \mu f_{4}{\rm e}^{-(\eta+1)} + \mu^{2} g_{4}{\rm e}^{-(\eta+1)},
\end{equation}
where
\begin{equation}\label{eq_rhs_J}
f_{4} =  (\xi^{2}-\eta)f_{2} - 2\xi f_{1},\quad
g_{4} = (\xi^{2}-\eta)g_{2} - 2\xi g_{1}.
\end{equation}
Then we introduce the Poincar\'e section $S_{-}$:
\begin{eqnarray}
S_{-} = \{(\,\xi, \eta, \zeta)\ : \ \xi =0, \ -1<\eta<0 ,\  \zeta\in \mathbb{R}\,\}.
\end{eqnarray}
We use the variables $(\zeta, J)$ as coordinates on this section and denote the first return map by 
\begin{equation}
P(\mu): S_{-}\to S_{-}.
\end{equation}
Since all trajectories of the unperturbed system are closed for $J>0$, 
the unperturbed Poincar\'e map $P(0)$ coincides with the identity map.
 It is natural to expect that after a perturbation a trajectory starting at $S_{-}$ will hit this section again near the initial point. 
Additionally one may expect that eigenvalues of the tangent map $TP(\mu)$ 
 should be close to one. However, this conjecture may become false if
 the trajectory is located near the separatrix where  the first return time grows substantially.
 To describe this effect we will concentrate our attention on the region
 near the separatrix where  $J\ll 1$.

Assuming $J_{0} = O(\mu )$, we consider a closed orbit  of the unperturbed system (\ref{eq_Unpert_v}), defined by $\zeta = \zeta_{0}$, $(\eta + 1 - \xi^2)e^{-(\eta+1)}=J_{0}$, and denote the orbit by $\mathcal{O}(\zeta_{0}, J_{0})$. We also introduce a 
$\mu \ln\mu-$neighborhood of the orbit, $U_{\mu}\left(\mathcal{O}(\zeta_{0}, J_{0})\right)$. Finally, we denote by $\mathcal{O}_{\mu}(\zeta_{0}, J_{0})$ the orbit of the system (\ref{eq_NormForm}), which contains the point $\mathcal{M}^{(0)}$ for which
\[ \xi =0,\quad \zeta =\zeta_{0},\quad J= J_{0} .\]
We suppose that $\mathcal{O}_{\mu}(\zeta_{0}, J_{0}) \subset U_{\mu}\left(\mathcal{O}(\zeta_{0}, J_{0})\right)$.

One may note that a trajectory corresponding to the unperturbed orbit $\mathcal{O}(\zeta_{0}, J_{0})$ possesses different behaviour in different regions.
It starts  at the point $\mathcal{M}^{(0)}$
and moves  initially near the separatrix $J=0$ till
it comes close to a turning point 
\begin{equation}
\xi_{+}=(k^{-1}-1)^{1/2},\quad \eta_{+}=k^{-1}-1,\quad \zeta_{+}=\zeta_{0},
\end{equation}
where 
\begin{equation} \label{eq_k}
k = \frac{1}{\ln {J_{0}}^{-1}},
\qquad k = O \left( \frac{1}{\ln \mu^{-1}} \right) \,.
\end{equation}
Then the trajectory "detaches" from the separatix, 
turns  in the $\xi$-direction and then "flies" across the region between two  branches of the separatrix. 
At the top of the orbit
\[ \xi_t =0, \quad \eta_t = k^{-1} +\ln k^{-1} -1 +o(1), \quad \zeta_t = \zeta_0.\]
Then the trajectory approaches 
the second turning point which is located symmetrically at $\xi_{-}=-(k^{-1}-1)^{1/2}$, $ \eta_{-}=k^{-1}-1$,  $\zeta_{-}=\zeta_{0}$, where the direction of motion with respect to 
$\xi$ is changed again  and finally the trajectory follows
the separtrix till an  intersection with $S_{-}$.  

According to this description we highlight in the neighborhood $U_{\mu}\left(\mathcal{O}(\zeta_{0}, J_{0})\right)$  the following overlapping domains:
\begin{align*}
&\mathcal{D}_{1} = \{(\xi, \eta, \zeta)\in U_{\mu}\left(\mathcal{O}(\zeta_{0}, J_{0})\right): \vert\xi\vert\ll \xi_{+},\; \eta<\eta_{+}\},\\
&\mathcal{D}_{2}^{\pm} = \{(\xi, \eta, \zeta)\in U_{\mu}\left(\mathcal{O}(\zeta_{0}, J_{0})\right): \vert\xi - \xi_{\pm}\vert\ll 1,\; \vert \eta - \eta_{\pm}\vert \ll k^{-1} \},\\
&\mathcal{D}_{3} = \{(\xi, \eta, \zeta)\in U_{\mu}\left(\mathcal{O}(\zeta_{0}, J_{0})\right): \vert\xi\vert\ll \xi_{+},\; \eta>\eta_{+}\}.
\end{align*}

Taking into account the symmetry of the system, we define an auxiliary Poincar\'e section
\begin{equation}\label{eq_S123}
S_{+}= \{(\xi, \eta, \zeta)\in U_{\mu}\left(\mathcal{O}(\zeta_{0}, J_{0})\right): \xi=0, \eta > \eta_{+} , \zeta\in \mathbb{R}\}
\end{equation}
and  introduce an auxiliary Poincar\'e map
\begin{eqnarray*}
F(\mu) : S_{-} \rightarrow S_{+}.
\end{eqnarray*}

\begin{figure}[h] 
\begin{center}
\includegraphics[width=8cm]{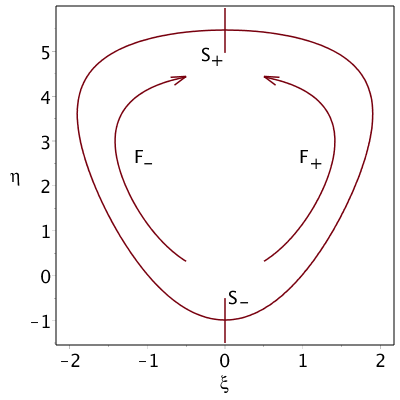}
\end{center}
\caption{Poincar\'e maps $F_{\pm} = F(\pm\mu)$.}
\label{Fig:maps_F_pm}
\end{figure}
We begin with considering the system (\ref{eq_NormForm}), (\ref{eq_dynJ}) in the domain $\mathcal{D}_{1}$ with initial conditions:
\[ \xi (0) = 0,\  \zeta (0) = \zeta_{0}, \  J(0)= J_{0},\ \eta (0) = \eta_0: \   \left( \eta_0 +1 \right) e^{-\left( \eta_0 +1 \right)} = J_{0}
\]
and we find some point $\mathcal{M}^{(1)}$ at the orbit
$\mathcal{O}_{\mu}(\zeta_{0}, J_{0})$ such that it belongs to 
$\mathcal{D}_{1} \cap \mathcal{D}_{2}^{+}$. Then we consider the system in $\mathcal{D}_{2}^{+}$ with corresponding initial conditions and find a point $\mathcal{M}^{(2)}$ at $\mathcal{O}_{\mu}(\zeta_{0}, J_{0})\cap \mathcal{D}_{2}^{+}\cap \mathcal{D}_{3}$. Finaly we consider the system in $\mathcal{D}_{3}$ and find the point $\mathcal{M}^{(3)}$ that belongs to $\mathcal{O}_{\mu}(\zeta_{0}, J_{0}) \cap S_{+}$. In this way we get $F(\mu )$. 

We also note that due to invariance of the system (\ref{eq_NormForm}) with respect to change (\ref{eq_Symm}) the map $F(-\mu)$ corresponds to the Poincar\'e map between sections 
$S_{-}$ and $S_{+}$, but backward in time (see fig.\ref{Fig:maps_F_pm}). 

Therefore the first-return map $P(\mu): S_{-} \rightarrow S_{-}$ can be represented as a composition
\begin{eqnarray}\label{eq_rep_P}
P(\mu) = F^{-1}(-\mu)\circ F(\mu ).
\end{eqnarray}

\subsection{Fixed points and the period-doubling bifurcation}

Since all trajectories of the unperturbed system are closed for $J>0$ one may expect that after a small perturbation the image of a point $(\zeta_0, J_0)\in S_{-}$ under the Poincar\'e map will be close to the initial point $(\zeta_0, J_0)$. The condition for a fixed point reads as
\begin{equation}\label{eq_fixed_point_1}
P(\mu )\left(
\begin{array}{cc}
\zeta_0 \\
J_0\end{array}
\right) 
= 
\left(
\begin{array}{cc}
\zeta_0 \\
J_0\end{array}
\right) 
\end{equation}
or, taking into account (\ref{eq_rep_P}), it can be rewritten
\begin{equation}\label{eq_fixed_point_2}
F(\mu)\left(
\begin{array}{cc}
\zeta_0 \\
J_0\end{array}
\right) 
= 
F(-\mu)\left(
\begin{array}{cc}
\zeta_0 \\
J_0\end{array}
\right).
\end{equation}
Due to smooth dependence of (\ref{eq_NormForm}) on $\mu$  one may represent the map $F(\mu)$ in the form
\begin{equation}\label{eq_rep_F}
F(\mu) = \mathrm{id} + \mu F_{1} + \mu^{2} F_{2} +\mu^3 F_3 + O(\mu^{4}).
\end{equation}
Then the symmetry of the normal form  (\ref{eq_fixed_point_2}) implies that 
\begin{equation}
\label{eq_fixed_point_3}
F_1, F_3: \left(
\begin{array}{cc}
\zeta_0 \\
J_0\end{array}
\right) 
\mapsto 
\left(
\begin{array}{cc}
0 \\
0\end{array}
\right).
\end{equation}
The period-doubling bifurcation occurs when one of the eigenvalues of the tangent map $D_{(\zeta_0, J_0)}P$ at a fixed point passes through $-1$. Thus, the condition for the period-doubling bifurcation can be written in the form
\begin{equation}
{\rm det}\big(D_{(\zeta_0, J_0)}P(\mu) + I\big) = 0
\end{equation}
or equivalently due to (\ref{eq_rep_P}) as
\begin{equation}\label{eq_period_doubling_1}
{\rm det}\big( D_{(\zeta_0, J_0)}F(\mu) + 
D_{(\zeta_0, J_0)}F(-\mu) \big) = 0.
\end{equation}
Substituting (\ref{eq_rep_F}) into (\ref{eq_period_doubling_1}), one gets
\begin{equation*}
{\rm det}\big( 2I + 2\mu^{2} D_{(\zeta_0, J_0)} F_2 + O(\mu^{4})\big) = 0
\end{equation*}
or
\begin{equation}\label{eq_period_doubling_2}
1+\mu^{2} {\rm Tr} \, D_{(\zeta_0, J_0)} F_2 + O(\mu^{4}) = 0.
\end{equation}
The condition (\ref{eq_period_doubling_2}) justifies the necessity of second order approximation 
(\ref{eq_rep_F}) for the map $F(\mu )$ to reveal the cascade of the period-doubling bifurcations.


\subsection{Domain $\mathcal{D}_{1}$}

In this subsection we consider the system (\ref{eq_NormForm}) in the domain $\mathcal{D}_{1}$.  As the orbit is close to the separatrix in this domain,
we introduce new variable $\vartheta $ instead of $\eta$: 
\begin{equation}\label{eq_Gamma}
\eta =\vartheta +\xi^2-1.
\end{equation}
Substituting (\ref{eq_Gamma}) into (\ref{eq_NormForm}) and (\ref{eq_dynJ}), one obtains
the system for $(\xi,\vartheta,\zeta ,J)$:
\begin{eqnarray}\label{eq_sys_D1}
\nonumber
\xi'&=&1-\vartheta
+\mu f_1 +
 \mu^2 g_1 \,, \\
\vartheta' &=&  2\xi\vartheta +\mu ( f_2 -2\xi f_1 )
 +   \mu^2 (g_2 -2\xi g_1)  \,, \\
\nonumber
 \zeta' &=&  \mu f_3 
 +   \mu^2 g_3  \,, \\
\nonumber
 J' &=&  \mu f_4 {\rm e}^{-\vartheta-\xi^{2}}
 +   \mu^2 g_4 {\rm e}^{-\vartheta-\xi^{2}}    \, ,
\end{eqnarray}
where all functions $f_{i}, g_{i}$ are defined in (\ref{eq_rhs_nf}), (\ref{eq_rhs_J}) and are evaluated at the point $(\xi, \vartheta + \xi^{2} -1, 
\zeta)$. 

One may conclude that as the variable $\vartheta \ll 1$ in $\mathcal{D}_{1}$, the derivative of the variable $\xi$ in (\ref{eq_sys_D1}) is positive in this domain. Taking this into account, we choose $\xi$ as a new independent variable and rewrite the system for other variables 
$\vartheta $, $\zeta$ and $J$ as functions of $\xi$:
\begin{eqnarray}
\nonumber
\frac{d\vartheta }{d \xi} &=& \frac{2\xi\vartheta +\mu ( f_2 -2\xi f_1 )
 +   \mu^2 (g_2 -2\xi g_1) }{1-\vartheta
+\mu f_1 +
 \mu^2 g_1 } 
 , \\ \label{Eq:systD1}
 \frac{d\zeta }{d \xi } &=& \frac{ \mu f_3 
 +   \mu^2 g_3}{1-\vartheta
+\mu f_1 +
 \mu^2 g_1 }   , \\ \nonumber
 \frac{dJ }{d \xi } &=& \frac{ \mu f_4 
 +   \mu^2 g_4}{1-\vartheta
+\mu f_1 +
 \mu^2 g_1 }\;{\rm e}^{-\vartheta -\xi^{2}}  .
\end{eqnarray}
We will find a solution of the system (\ref{Eq:systD1}) as a perturbation of the orbit $\mathcal{O}(\zeta_0, J_0)$
\begin{eqnarray}
\nonumber
 \vartheta (\xi )&=& \vartheta_0(\xi )  +\mu \vartheta_1 (\xi ) +\mu^2\vartheta_2 (\xi ) +O(\mu^3), \\ 
\label{eq_Expand_1}
\zeta (\xi )  &=& \zeta_0 +\mu \zeta_1(\xi )  +\mu^2\zeta_2 (\xi ) +O(\mu^3),\\ 
\nonumber
J(\xi ) &=& J_0+\mu J_1 (\xi ) +\mu^2 J_2(\xi )  +O(\mu^3)  
\end{eqnarray}
and fix initial conditions for the unknown functions in the following way:
\begin{equation}\label{eq_InitCond_1bb}
\vartheta_{i} (0)=0,\quad 
\zeta_{i}(0)=0,\quad
J_{i}(0)=0,\quad i=1,2.
\end{equation}

The function $\vartheta_0 (\xi )$ corresponds to the unperturbed orbit and is given implicitly by the following equation
\begin{equation}\label{eq_gamma_0}
\vartheta_0 e^{-\vartheta_0}=J_{0} e^{\xi^2}.
\end{equation}
Hence
\begin{equation}\label{eq_est_gamma_0}
\vartheta_0(\xi )=J_{0} e^{\xi^2} \left( 1 +O \left( J_{0} e^{\xi^2}
\right) 
\right).
\end{equation}

Now we fix the point $\mathcal{M}^{(1)}$ by setting the corresponding value of $\xi =\xi^{(1)}$ as
\begin{equation}\label{eq_Xi1}
\xi^{(1)} = k^{-1/2}\left(1 - \frac{1}{2}k\ln k^{-1}\right) .
\end{equation}
Then, taking into account (\ref{eq_k}), one obtains that 
\begin{equation}\label{eq_Estim_S_1}
(\xi^{(1)})^2 = k^{-1}-\ln k^{-1} +\frac14 k\ln^2 k^{-1},\qquad
e^{(\xi^{(1)})^2} =\frac{k}{J_{0}} \left(1+O(k\ln^2 k^{-1} )\right).
\end{equation} 
and the following  estimates hold  
\begin{equation}\label{eq_est1_gamma_0}
\vartheta_0 (\xi ) = J_0 e^{\xi^2} +O(k^2), \qquad   \vartheta_{0} (\xi ) = O(k).
\end{equation}

Substituting (\ref{eq_Expand_1}) into (\ref{Eq:systD1}) and expanding with respect to 
$\mu$, one obtains equations for the functions $\vartheta_{i}, \zeta_{i}, J_{i}$. Note that
derivatives ${\rm d}\zeta/{\rm d}\xi$, ${\rm d}J/{\rm d}\xi$ are of the order $O(\mu)$. This implies the function $\vartheta_{2}$ does not appear in equations for $\zeta_{2}$, $J_{2}$ and we need to construct the solution $\vartheta (\xi)$ only up to terms of the order $O(\mu)$.

The function $\vartheta_1 (\xi)$ satisfies an equation:
\[ 
\frac{d\vartheta_{1}}{d \xi} = \frac{2\xi \vartheta_1}{(1-\vartheta_0)^2 } 
+ \mathcal{R}_{1}^{(\vartheta )}(\xi),\]
where
\[ 
\mathcal{R}_{1}^{(\vartheta )}(\xi) =
\frac{(1-\vartheta_{0})f_{2} - 2\xi f_{1}}{(1-\vartheta_0)^2 },\]
and the functions $f_{1}, f_{2}$ are evaluated at the point 
$(\xi, \vartheta_{0} +\xi^{2} -1,  \zeta_0) $.

Taking into account (\ref{eq_InitCond_1bb}), the solution of the equation on $\vartheta_1(\xi )$ has the form
\[
\vartheta_1 (\xi ) = e^{\int_0^\xi \frac{2s\,ds}{(1-\vartheta_0(s))^2}} 
 \int_0^\xi  e^{-\int_0^s \frac{2p\,dp}{(1-\vartheta_0(p))^2}}\cdot
\mathcal{R}_{1}^{(\vartheta )}(s)ds  .\]

From (\ref{eq_est1_gamma_0}) we get
\[ \int_0^s \frac{2p\,dp}{(1-\vartheta_0(p))^2}
= \int_0^s 2p\big( 1 +O(\vartheta_0(p))\big) dp
 = s^2 + O(k).
\]
Then, using (\ref{eq_rhs_nf}) and (\ref{eq_est1_gamma_0}), one obtains
\[ \vartheta_1 (\xi ) = e^{\xi^2}  \int_0^\xi e^{-s^2}  \left[ -2\gamma s^4 +(\alpha_1 -2\gamma_0 )s^2  - \alpha_1+  \alpha_2 \zeta_0 
 \right]{\rm d}s 
\cdot\left( 1+O(k) \right).
\]
Note that 
\[ \Phi (\xi) =  \int_0^\xi e^{-s^2} ds= \frac{\sqrt{\pi}}{2} +O \left(
\frac{e^{-\xi^2}}{\xi }\right) , \quad \xi\to +\infty , \]
\[  \int_0^\xi s^2 e^{-s^2} ds= -\frac12 \xi e^{-\xi^2} +\frac12  \Phi (\xi) ,\]
\[ \int_0^\xi s^4 e^{-s^2} ds= -\frac12 \xi^3 e^{-\xi^2}
-\frac34 \xi e^{-\xi^2}  +\frac34  \Phi (\xi) .\]
We introduce  
\begin{equation}\label{eq_def_A}
A(\zeta_{0}) = \frac{\sqrt{\pi }}{2}
\left( - \frac32 \gamma + \alpha_2\zeta_{0} -\frac12 \alpha_1 - \gamma_0  \right) .
\end{equation}
Thus, for $\vartheta_1 (\xi )$ one has
\begin{equation}
  \label{eq_vartheta_1}
  \vartheta_1 (\xi ) =\left[ e^{\xi^2}  \frac{2}{\sqrt{\pi}} A(\zeta_{0}) \Phi (\xi) 
+Q_3^{(\vartheta )}(\xi) \right] (1+O(k)),
\end{equation} 
where $Q_3^{(\vartheta )} (\xi )$ is a polynomial of the third order in $\xi$.
Using (\ref{eq_Estim_S_1}), we may calculate $\vartheta_1 (\xi)$ at the point 
$\mathcal{M}^{(1)}$
\begin{equation} \label{eq_theta1_S1}
\vartheta_1^{(1)}= \vartheta_1 (\xi^{(1)}) = A(\zeta_{0}) 
\frac{k}{J_{0}} \Bigl(1+O(k \ln^2 k^{-1})\Bigr) +O(k^{-3/2}) .
\end{equation}

From (\ref{Eq:systD1}) one gets that the $i$-th order approximation $(\zeta_{i}, J_{i})$ satisfies a system of equations  ($i=1,2$):
\begin{eqnarray}
\label{Eq:syst_k_D1}
 \frac{d\zeta_{i} }{d \xi } &=& \mathcal{R}_{i}^{(\zeta)}(\xi), \\ 
\nonumber
\frac{d J_{i}}{d \xi } &=& \mathcal{R}_{i}^{(J)}(\xi).
\end{eqnarray}
Here for $i=1$  
\begin{eqnarray}
\label{eq_rhs_1_D1}
\mathcal{R}_{1}^{(\zeta)}(\xi) &=& 
\frac{f_{3}}{1-\vartheta_0},\\
\nonumber
\mathcal{R}_{1}^{(J)}(\xi) &=& 
\frac{f_{4}}{1-\vartheta_0}{\rm e}^{-\vartheta_{0}-\xi^{2}} 
\end{eqnarray}
and for $i=2$
\begin{eqnarray}
\nonumber
\mathcal{R}_{2}^{(\zeta)}(\xi) &=& 
\frac{f'_{3,\eta}\vartheta_{1} + f'_{3,\zeta}\zeta_{1} +  g_{3}}{1-\vartheta_0} + 
\frac{f_{3}(\vartheta_{1} - f_{1})}{(1-\vartheta_0)^{2}},\\
\label{eq_rhs_2_D1}
\mathcal{R}_{2}^{(J)}(\xi) &=& \left [
\frac{f'_{4,\eta}\vartheta_{1} + f'_{4,\zeta}\zeta_{1} +  g_{4}}{1-\vartheta_0} + 
\frac{f_{4}(\vartheta_{0}\vartheta_{1} - f_{1})}{(1-\vartheta_0)^{2}}\right ]
{\rm e}^{-\vartheta_{0}-\xi^{2}} 
\end{eqnarray}
and the functions $f_{i}, g_{i}$, $i=1,\ldots, 4$ are evaluated at the point  
$(\xi, \vartheta_{0} +\xi^{2} -1, \zeta_0)$.

Taking into account (\ref{eq_InitCond_1bb}), the solution of the system (\ref{Eq:syst_k_D1}) has the form
\begin{eqnarray}
\label{eq_sol_k_D1}
\zeta_{i}(\xi) &=&  \int_{0}^{\xi}\mathcal{R}_{i}^{(\zeta)}(s){\rm d}s,\\
\nonumber
J_{i}(\xi) &=&  \int_{0}^{\xi}\mathcal{R}_{i}^{(J)}(s){\rm d}s.
\end{eqnarray}

Substituting (\ref{eq_rhs_nf}) into (\ref{eq_rhs_1_D1}) and then into (\ref{eq_sol_k_D1}) and taking into account (\ref{eq_est1_gamma_0}), we get
\begin{eqnarray}
\nonumber
\zeta_1(\xi) &=&   \int_0^\xi \left[\beta_{1}(s^2-1) + \beta_2 \zeta_0    \right]
{\rm d}s  \cdot\left( 1+O(k) \right),\\
\nonumber
J_{1}(\xi) &=&   \int_0^\xi {\rm e}^{-s^2}\left[ - 2 \gamma s^4 +(\alpha_1 -2\gamma_0 )s^2 -
\alpha_1 +  \alpha_2 \zeta_0 \right]{\rm d}s 
\cdot\left( 1+O(k) \right).
\end{eqnarray}
Consequently, 
\begin{eqnarray}
\nonumber
\zeta_1(\xi) &=& \left[ \frac13 \beta_1 \xi^3- \beta_{1} \xi
 + \beta_2 \zeta_{0} \xi \right] (1 +O(k)) ,\\
\label{eq_as_1_D1}
J_1 (\xi) &=& \left[ \frac{2}{\sqrt{\pi}}
 A(\zeta_{0}) \Phi (\xi)  + e^{-\xi^2} Q_3^{(J )}(\xi) \right]  \Bigl(1+O(k )\Bigr).
\end{eqnarray}  
where $Q_3^{(J )} (\xi )$ is a polynomial of the third order in $\xi$.

Then, using (\ref{eq_Xi1}), we have at the point $\mathcal{M}^{(1)}$
\begin{eqnarray}
\nonumber
\zeta_1^{(1)}=\zeta_1(\xi^{(1)}) &=&  k^{-3/2} \left[ \frac13 \beta_1 + \beta_2 k\zeta_{0} - \frac{1}{2}\left(\beta_{1}+\beta_{2}k\zeta_{0}\right)k \ln k^{-1} - \beta_{1}k + O(k^{2} \ln^{2} k^{-1})\right],\\
\label{eq_as_1_S1}
J_1^{(1)}=J_1 (\xi^{(1)}) &=&  A(\zeta_{0})\Bigl(1+O(k )\Bigr) +O \left( \frac{J_0}{k^{5/2}}  \right).
\end{eqnarray}

For $i=2$ we use (\ref{eq_rhs_2_D1}) and
take into account (\ref{eq_vartheta_1}), (\ref{eq_as_1_D1}) to obtain
\[ \mathcal{R}_{2}^{(\zeta)}(\xi) =\left[ (\beta_1 \xi^2  + \beta_2 \zeta_0 )  e^{\xi^2} \frac{2}{\sqrt{\pi}} A(\zeta_{0}) \Phi (\xi)  + Q_5^{(\zeta )} (\xi )
\right] (1+O(k)),
\]
where $Q_5^{(\zeta )} (\xi )$ is a polynomial of the fifth order in $\xi$.
Note that due to (\ref{eq_Estim_S_1})
\[   \int_0^{\xi^{(1)}} e^{s^2} \Phi (s) \, ds = \frac{\sqrt{\pi}}{4} \frac{k^{3/2}}{J_0} (1+O(k \ln^2 k^{-1})) , \]
\[  \int_0^{\xi^{(1)}} s^2 e^{s^2} \Phi (s) \, ds = \frac{\sqrt{\pi}}{4} \frac{k^{1/2}}{J_0} (1+O(k \ln^2 k^{-1})) .
\] 
Then for $i=2$ we get from (\ref{eq_sol_k_D1}) 
\begin{equation}
  \label{eq_as_zeta2_S1}
\zeta_2^{(1)} = \zeta_{2}(\xi^{(1)}) =  A(\zeta_{0})  
\frac{k^{1/2}}{2J_{0}} \left[  \beta_1 + \beta_2 k\zeta_{0} \right] 
\Bigl(1+O(k \ln^2 k^{-1})\Bigr) + O(k^{-3}) .
\end{equation}

Application of formulae (\ref{eq_rhs_2_D1}), (\ref{eq_as_1_D1}) and (\ref{eq_vartheta_1}) yields 
\[ \mathcal{R}_{2}^{(J)}(\xi) =\left[ 
e^{-\xi^2} Q_7^{(J)} (\xi ) 
- ( \alpha_1 \xi^2 + \alpha_2 \zeta_0 - 2\alpha_1
)
 \frac{2}{\sqrt{\pi}} A(\zeta_0) \Phi (\xi )
\right] (1+O(k)),
\]
where $Q_7^{(J )} (\xi )$ is a polynomial of the seventh order in $\xi$.
Note that due to (\ref{eq_Xi1})

\[ \int_0^{\xi^{(1)}} \Phi (s) ds = \frac{\sqrt{\pi}}{2} k^{-1/2} (1+O(k\ln k^{-1})) - \frac12 ,
\]
\[ \int_0^{\xi^{(1)}} s^2 \Phi (s) ds = \frac{k^{-3/2}}{3} \frac{\sqrt{\pi}}{2} (1+O(k\ln k^{-1})) .\] 

Then we get 
\begin{equation}\label{eq_as_J2_S1}
 J_2^{(1)}= J_{2} (\xi^{(1)}) =  
- k^{-3/2} A(\zeta_{0}) 
 \left( \frac{\alpha_{1}}{3} +\alpha_2 k \zeta_0 \right) \Bigl(1+O(k\ln k^{-1} )\Bigr) +O(k^{-1}).
\end{equation}

Thus, the coordinates $(\zeta^{(1)}, J^{(1)})$ of the point $\mathcal{M}^{(1)}$ which belongs to $\mathcal{O}_{\mu}(\zeta_{0}, J_{0}) \cap \mathcal{D}_{1} \cap \mathcal{D}_{2}^{+}$ are described by
\begin{align}\label{eq_as_S1}
&\zeta^{(1)} = \zeta_{0} + \mu \zeta_{1}^{(1)} + \mu^{2} \zeta_{2}^{(1)}+O(\mu^3),\\
\nonumber
&J^{(1)} = J_{0} + \mu J_{1}^{(1)} + \mu^{2} J_{2}^{(1)}+O(\mu^3), 
\end{align}
where $\zeta_i^{(1)}$ and $J_i^{(1)}$ are defined by
(\ref{eq_as_1_S1}), (\ref{eq_as_zeta2_S1}) and (\ref{eq_as_J2_S1}).


\subsection{Domain $\mathcal{D}_{2}^{+}$}

The aim of this subsection is to obtain the second order approximation for a point $\mathcal{M}^{(2)}=(\zeta^{(2)}, J^{(2)}) \in \mathcal{O}_{\mu}(\zeta_{0}, J_{0}) \cap \mathcal{D}_{2}^{+} \cap \mathcal{D}_{3}$.
We
fix the point $\mathcal{M}^{(2)}$ by setting a value of the variable $\eta$ corresponding to this point as 
\begin{equation}
\eta^{(2)} = k^{-1} + \frac{1}{2}\ln k^{-1} - 1 .
\end{equation}
To consider the system in the domain $\mathcal{D}_{2}^{+}$ it is convinient to introduce new variables $z$ and $v$ by the following way:
\[ \xi = k^{-1/2}(1-kz),\qquad \eta = k^{-1}+ v  -1 . \]
Then the coordinate $v$ corresponding to $\mathcal{M}^{(2)}=(\zeta^{(2)}, J^{(2)})$ is
\begin{equation}
  \label{eq_S2}
  v^{(2)}=\frac12 \ln k^{-1} . 
\end{equation}
On the other hand, the coordinates $(z,v )$ which correspond to the point $\mathcal{M}^{(1)}=(\zeta^{(1)}, J^{(1)})$ are
\[ z^{(1)}=\frac12 \ln k^{-1} , \quad 
v^{(1)} = \eta^{(1)}-k^{-1}+1=\vartheta^{(1)}+(\xi^{(1)})^2-k^{-1} . \]
Therefore the coordinate $v^{(1)}$ depends on $\mu$:
\[ v^{(1)} = v_0^{(1)} + \mu v_1^{(1)} + O(\mu^2),\]
where, due to relations (\ref{eq_Estim_S_1}), (\ref{eq_est1_gamma_0}),  (\ref{eq_theta1_S1}), 
\begin{equation} \label{Eq:v01}
v_{0}^{(1)} = -\ln k^{-1} + O(k\ln^{2}k^{-1}),
\end{equation}
\begin{eqnarray}
\label{eq_v1}
v_{1}^{(1)} = \vartheta_{1}^{(1)} =  
 A(\zeta_{0}) \frac{k}{J_{0}}\cdot \Bigl(1 + O(k\ln^2 k^{-1})\Bigr)
 +O(k^{-3/2}).
\end{eqnarray}

In terms of new variables the equations of motion (\ref{eq_NormForm}) and (\ref{eq_dynJ}) can be rewritten as
\begin{eqnarray}
\nonumber
z'&=& 2k^{-1/2}z-k^{1/2}z^2+k^{-1/2}(v-1)-\mu k^{-1/2}f_1 - \mu^2 k^{-1/2} g_1  \,, \\
\nonumber
v'&=& 2k^{-1/2} (1-kz)+\mu f_2+\mu^2g_2 \,, \\
\nonumber
 \zeta ' &=&  \mu f_3
 +   \mu^2 g_3   \, ,\\
 \nonumber
J ' &=&  J_{0}\left[\mu f_4
 +   \mu^2 g_4\right]  {\rm e}^{-v}\, ,
\end{eqnarray}
where the functions $f_{i}, g_{i}, i=1,\ldots, 4$ are evaluated at the point 
$(k^{-1/2}(1-kz), k^{-1}-1+v, \zeta)$.

Note that in the domain $\mathcal{D}_{2}^{+}$ an expression $kz< 1$ and the derivative $v'$ does not vanish. Thus, one can take $v$ as a new independent variable and obtain equations on $(z, \zeta, J)$ as functions of $v$:
\begin{eqnarray}
\nonumber
\frac{dz}{dv}&=& \frac{2z+v-1-kz^2-\mu f_1-\mu^2 g_1}{2(1-kz)+\mu k^{1/2}f_2+\mu^2 k^{1/2} g_2} \,, \\
\label{eq_D2}
\frac{d\zeta}{dv} &=& \frac{k^{1/2}(\mu f_3+ \mu^2 g_3)}{2(1-kz)+\mu k^{1/2}f_2+\mu^2 k^{1/2} g_2}  \,, \\
\nonumber
\frac{d J}{dv} &=& \frac{J_0 k^{1/2}(\mu f_4+ \mu^2 g_4)}{2(1-kz)+\mu k^{1/2}f_2+\mu^2 k^{1/2} g_2]}{\rm e}^{-v}  \, .
\end{eqnarray}
We supply (\ref{eq_D2}) by initial conditions corresponding to the point $\mathcal{M}^{(1)}$:
\begin{equation}\label{eq_InitCond_2}
z(v^{(1)}) = z^{(1)},\quad
\zeta(v^{(1)}) = \zeta^{(1)},\quad
J(v^{(1)}) = J^{(1)}
\end{equation}
and  find the solution satisfying (\ref{eq_D2}) and (\ref{eq_InitCond_2}) in a form
\begin{eqnarray}
\nonumber
z(v)&=& z_0(v)  +\mu z_1 (v)  +O(\mu^2), \\ 
\label{eq_Expand_2}
\zeta (v)  &=& \zeta_{0} +\mu \zeta_1(v)  +\mu^2\zeta_2 (v) +O(\mu^3),\\ 
\nonumber
J(v) &=& J_{0}+\mu J_1 (v) +\mu^2 J_2(v)  +O(\mu^3) , 
\end{eqnarray}
where $z_0(v)$ corresponds to the unperturbed orbit $\mathcal{O}(\zeta_{0}, J_{0})$.

The initial conditions (\ref{eq_InitCond_2}) are set at the point $v^{(1)}$ which depends on $\mu$.
Using the Taylor formula with respect to $\mu$, one may reformulate these conditions at the point $v_0^{(1)}$ as follows:
\begin{align}\label{eq_InitCond_2b}
\nonumber
 & z_0(v_0^{(1)}) = z_{0}^{(1)} ,\quad z_{1}(v_0^{(1)}) =  - \frac{d z_{0}}{d v}(v_{0}^{(1)})\cdot v_{1}^{(1)},\\
& \zeta_{1}(v_0^{(1)})  = \zeta_{1}^{(1)},\quad
\zeta_{2}(v_0^{(1)})  = \zeta_{2}^{(1)} - \frac{d \zeta_{1}}{d v}(v_{0}^{(1)})\cdot v_{1}^{(1)},\\
\nonumber
& J_{1}(v_0^{(1)}) =  J_{1}^{(1)},\quad J_{2}(v_0^{(1)})  = J_{2}^{(1)} - \frac{d J_{1}}{d v}(v_{0}^{(1)})\cdot v_{1}^{(1)} .
\end{align} 

Substituting (\ref{eq_Expand_2}) into (\ref{eq_D2}) and collecting the terms of the same order of $\mu$, one gets equations for the components $z_{i}, \zeta_{i}, J_{i}$. 
We solve these equations with the initial conditions (\ref{eq_InitCond_2b}) to obtain 
$\zeta_i(v)$ and $J_i(v)$. Finally, substituting $v=v^{(2)}$, one gets the point $\mathcal{M}^{(2)}  =(\zeta^{(2)}, J^{(2)})$.

Thus, our task is to derive asymptotic formulae for $\zeta_{1,2} (v^{(2)})$ and $J_{1,2}(v^{(2)})$. We begin with auxilary asymptotics for $z_0(v)$ and $z_1(v)$.

Note that $z_0(v)$ corresponds to the unperturbed orbit $\mathcal{O}(\zeta_0, J_0)$. Hence, due to (\ref{eq_unpertOrbit}), it is a solution of the following equation
\begin{equation}\label{eq_z0_sol}
z_{0} - \frac{1}{2} k z_{0}^{2} = \frac{1}{2}\Bigl({\rm e}^{v} - v\Bigr).
\end{equation}

Since $k ({\rm e}^{v} - v)\ll 1$ in the domain $\mathcal{D}_{2}^{+}$, the function $z_0(v)$ admits the following asymptotics
\begin{equation}\label{eq_asymp_z0_D2}
z_{0} (v) = \frac{1}{2}\Bigl({\rm e}^{v} - v \Bigr)
\Bigl(
1  + O\left( k\left({\rm e}^{v} - v \right) \right) \Bigr).
\end{equation}

It is not difficult to verify that in $\mathcal{D}_2^{+}$
\begin{equation} \label{Eq:asforz0}
  k z_0(v) =O(k^{1/2}).
\end{equation}
 
From (\ref{eq_D2}) one may deduce that equation for $z_1$ as a function of $v$ can be written as
\[ 
\frac{d z_{1}}{d v} = \left(1 + k\frac{dz_{0}/dv}{1-kz_{0}}\right) z_{1} 
+ \mathcal{R}_{1}^{(z)}(v),
\]
where
\[ \mathcal{R}_{1}^{(z)}(v)
= - \frac{f_{1}}{2(1-k z_{0}(v))} - 
k^{1/2}\frac{(2z_{0}(v) + v- 1 -kz_0^2) f_{2}}{4(1-k z_{0}(v))^{2}}
 = O( k^{-3/2} ).
\]
Then
\begin{equation}\label{Eq:z1sol}
z_1 (v) = e^{\int_{v_0^{(1)}}^{v}\left [1 + k\frac{dz_{0}/ds}{1-kz_{0}} \right]ds}
\left( z_{1}\left( v_0^{(1)}\right)
 + \int_{v_0^{(1)}}^{v}  
e^{-\int_{v_0^{(1)}}^{s}\left [1 + k\frac{dz_{0}/dp}{1-kz_{0}} \right]dp}\cdot
\mathcal{R}_{1}^{(z)}(s)ds \right) .  
\end{equation}
Due to (\ref{Eq:v01}) and (\ref{Eq:asforz0}) we have
\[ \int_{v_0^{(1)}}^{v}\left [1 + k\frac{dz_{0}/ds}{1-kz_{0}} \right]ds = v+ \ln k^{-1} +O(k^{1/2}).
\]
Then, using  (\ref{eq_InitCond_2b}), (\ref{eq_asymp_z0_D2})  and (\ref{eq_v1}), one gets
\begin{equation}\label{eq_z1}
  z_1(v)=  \frac{e^v}{2J_0}A(\zeta_0)
 \, (1+O(k^{1/2})).
\end{equation}
Finally, (\ref{eq_S2}) yields
\begin{equation}
  \label{eq_z1_S2}
   z_1(v^{(2)}) = \frac{k^{-1/2}}{2J_0} A(\zeta_0)
(1+O(k^{1/2})).
\end{equation}

We substitute (\ref{eq_Expand_2}) into (\ref{eq_D2}) and obtain equations for $\zeta_i$ and $J_i$ ($i=1,2$) in the following form
\begin{eqnarray}
\label{Eq:syst_m_D2}
\frac{d\zeta_{i} }{d v } &=& \mathcal{R}_{i}^{(\zeta)}(v), \\ 
\nonumber
\frac{d J_{i}}{d v } &=& \mathcal{R}_{i}^{(J)}(v),
\end{eqnarray}
where for $i=1$:
\begin{eqnarray}
\label{eq_rhs_1_D2}
\mathcal{R}_{1}^{(\zeta)}(v) &=& 
k^{1/2}\frac{f_{3}}{2(1-k z_{0}(v))}
,\\
\nonumber
\mathcal{R}_{1}^{(J)}(v) &=& 
J_0k^{1/2}\frac{f_{4}}{2(1-k z_{0}(v))}{\rm e}^{-v} ,
\end{eqnarray}
and for $i=2$:
\begin{eqnarray}
\label{eq_rhs_2_D2}
\mathcal{R}_{2}^{(\zeta)}(v) &=& 
k^{1/2}\left[ 
\frac{ - k^{1/2} f'_{3,\xi} \cdot z_{1}(v) + f'_{3,\zeta}\cdot \zeta_{1}(v) +  
g_{3}}{2(1-k z_{0}(v))} + 
\frac{f_{3}\cdot (2 k z_{1}(v) - k^{1/2} f_{2})}{4(1-k z_{0}(v))^{2}}\right],\\
\nonumber
\mathcal{R}_{2}^{(J)}(v) &=& 
J_{0} k^{1/2}\left[ 
\frac{ - k^{1/2} f'_{4,\xi} \cdot z_{1}(v) + f'_{4,\zeta}\cdot \zeta_{1}(v) +  
g_{4}}{2(1-k z_{0}(v))} + 
\frac{f_{4}(2 k z_{1}(v) - k^{1/2}f_{2})}{4(1-k z_{0}(v))^{2}}\right]
{\rm e}^{-v}.
\end{eqnarray}
Here the functions $f_{i}, g_{i}$ are evaluated at the point $(\xi, \eta, \zeta) = \Bigl(k^{-1/2}(1-k z_{0}(v)), k^{-1} + v - 1, \zeta_{0}\Bigr)$.

The solutions of (\ref{Eq:syst_m_D2}) are
\begin{eqnarray}
\label{eq_sol_m_D2}
\zeta_{i}(v) &=& \zeta_{i}\left( v_0^{(1)}\right) +\int_{v_0^{(1)}}^{v}
\mathcal{R}_{i}^{(\zeta)}(s){\rm d}s,\\
\nonumber
J_{i}(v) &=&  J_{i}\left( v_0^{(1)}\right)+\int_{v_0^{(1)}}^{v}\mathcal{R}_{i}^{(J)}(s){\rm d}s.
\end{eqnarray}

Due to (\ref{eq_rhs_nf}), (\ref{eq_rhs_J}) and the estimate (\ref{Eq:asforz0}) one may deduce from (\ref{eq_rhs_1_D2}) that the functions $ \mathcal{R}_{1}^{(\zeta)}(v)$ and $\mathcal{R}_{1}^{(J)}(v) $  can be written as
\begin{eqnarray}
\nonumber
\mathcal{R}_{1}^{(\zeta)}(v) 
&=& \frac12 k^{-1/2} (\beta_1 + \beta_2 k \zeta_0) + O(k^{1/2}\ln k^{-1}) ,\\ 
\label{eq_as_R_1_D2}
\mathcal{R}_{1}^{(J)}(v) 
&=& -\gamma J_0 k^{-3/2} e^{-v} (1+O(k^{1/2}))  
.
\end{eqnarray}

Substituting $v=v^{(2)} = \frac{1}{2}\ln k^{-1}$  into (\ref{eq_sol_m_D2})
 and taking into account the formulae
(\ref{eq_as_R_1_D2}), (\ref{eq_InitCond_2b}), (\ref{eq_as_1_S1}), one obtains
\begin{align}
\nonumber
&\zeta_1^{(2)}= \zeta_{1}(v^{(2)}) =  k^{-3/2} \left[ \frac13 \beta_1 + \beta_2 k\zeta_{0} + \frac{1}{4}\left(\beta_{1}+\beta_{2}k\zeta_{0}\right)k \ln k^{-1} - \beta_1 k +
O(k^{3/2}) \right]
,\\
\label{eq_as_1_S2}
&J_1^{(2)}=J_{1}(v^{(2)}) =  A(\zeta_{0})\Bigl(1+O(k )\Bigr)
 + O(J_0 k^{-5/2} ).
\end{align}

The formula (\ref{eq_sol_m_D2}) for $i=2$ together with (\ref{eq_InitCond_2b}) leads to
\begin{eqnarray}
\label{eq_bc22}
\zeta_{2}(v^{(2)}) &=& \zeta_2^{(1)} - \frac{d \zeta_{1}}{d v}(v_{0}^{(1)})\cdot v_{1}^{(1)} +\int_{v_0^{(1)}}^{v^{(2)}}
\mathcal{R}_{2}^{(\zeta)}(s){\rm d}s,\\
\nonumber
J_{2}(v^{(2)}) &=&  J_2^{(1)} - \frac{d J_{1}}{d v}(v_{0}^{(1)})\cdot v_{1}^{(1)} +\int_{v_0^{(1)}}^{v^{(2)}}
\mathcal{R}_{2}^{(J)}(s){\rm d}s.
\end{eqnarray}
Taking into account (\ref{eq_v1}), (\ref{Eq:syst_m_D2}) and (\ref{eq_as_R_1_D2}),
 one gets 
 \begin{eqnarray}
\nonumber
 - \frac{d \zeta_{1}}{d v}(v_{0}^{(1)})\cdot v_{1}^{(1)} &=&  
- A(\zeta_0)\frac{k^{1/2}}{2J_0}(\beta_1 +\beta_2k\zeta_0)(1+O(k^{1/2})) ,\\
\nonumber
- \frac{d J_{1}}{d v}(v_{0}^{(1)})\cdot v_{1}^{(1)} &=&  \gamma  k^{-3/2} A(\zeta_0) (1+O(k^{1/2})) .
\end{eqnarray}
Then (\ref{eq_rhs_2_D2}),  (\ref{Eq:asforz0}), (\ref{eq_z1})
yield
\begin{eqnarray}
\nonumber
\mathcal{R}_{2}^{(\zeta)}(v) 
&=& \frac{ k^{1/2}}{4 J_0 } A(\zeta_0)
(\beta_1 +\beta_2 k\zeta_0) e^v \Bigl(1+O(k^{1/2})\Bigr)
,\\ 
\nonumber
\mathcal{R}_{2}^{(J)}(v) 
&=& O( k^{-1/2} ) .
\end{eqnarray}
Hence
\begin{eqnarray}
\nonumber
\int_{v_0^{(1)}}^{v^{(2)}}\mathcal{R}_{2}^{(\zeta)}(s)\, ds 
&=& \frac{ 1}{4 J_0 } A(\zeta_0)
(\beta_1 +\beta_2 k\zeta_0)  \Bigl(1+O(k^{1/2})\Bigr)
,\\ 
\nonumber
\int_{v_0^{(1)}}^{v^{(2)}} \mathcal{R}_{2}^{(J)}(s) \, ds
&=& O( k^{-1/2} \ln k^{-1} ) .
\end{eqnarray}

Substituting these formulae together with (\ref{eq_as_zeta2_S1}), (\ref{eq_as_J2_S1}) into (\ref{eq_bc22}), we get
\begin{eqnarray}
\label{eq_as_2_S2}
\zeta_2^{(2)} =\zeta_{2}(v^{(2)}) &=&  A(\zeta_{0})\frac{1}{4J_0}\left[  \beta_1 + \beta_2 k\zeta_{0} \right] 
\Bigl(1+O(k^{1/2})\Bigr)
,\\
\nonumber
J_2^{(2)} =J_{2} (v^{(2)}) &=&  -A(\zeta_{0})k^{-3/2} \left( \frac{\alpha_{1} }{3} +\alpha_2 k \zeta_0 - \gamma \right)\Bigl(1+O(k^{1/2})\Bigr).
\end{eqnarray}

Consequently, the point $\mathcal{M}^{(2)}$ is characterized by
\begin{eqnarray}
  \label{eq_M2}
  \zeta^{(2)}=\zeta_0 +\mu \zeta_1^{(2)}
+\mu^2 \zeta_2^{(2)} +O(\mu^3), \\ \nonumber
J^{(2)}=J_0 +\mu J_1^{(2)}
+\mu^2 J_2^{(2)} +O(\mu^3),
\end{eqnarray}
where $\zeta_i^{(2)}$ and $J_i^{(2)}$ are defined by (\ref{eq_as_1_S2}) for $i=1$ and (\ref{eq_as_2_S2}) for $i=2$.


\subsection{Domain $\mathcal{D}_{3}$}

In this subsection we derive an asymptotic for the Poincar\'e map $F(\mu)$. In the domain $\mathcal{D}_{3}$ it is convenient to introduce variables
\begin{equation}\label{eq_var_D3}
y = k^{1/2}\xi,\qquad v = \eta - k^{-1} + 1.
\end{equation}
In terms of the variables $(y, v, \zeta)$ the section $S_{+}$ can be rewritten as
\begin{equation}
S_{+} = \{(y, v, \zeta)\in U_{\mu}\left(\mathcal{O}(\zeta_{0}, J_{0})\right): y=0,\; v > 0\}
\end{equation}
and the equations of motion (\ref{eq_NormForm}) and (\ref{eq_dynJ}) take the form
\begin{eqnarray}
\nonumber
y'&=& k^{-1/2} (y^2 -1 - k(v-1)) +\mu k^{1/2}f_1    +\mu^2 k^{1/2}  g_1   \,, \\
\nonumber
v'&=&   2k^{-1/2}y +\mu f_2  +\mu^2  g_2   \,, \\
\nonumber
 \zeta ' &=&  \mu f_3  
 +   \mu^2 g_3    \,, \\
\nonumber
 J ' &=&  J_0( \mu f_4  
 +   \mu^2 g_4 ) \, {\rm e}^{-v} \, ,
\end{eqnarray}
where the functions $f_{i}, g_{i}$ are evaluated at the point 
$(\xi, \eta, \zeta) = (k^{-1/2}y, v + k^{-1} - 1, \zeta)$.
One may note that in the domain $\mathcal{D}_{3}$ the derivative $y'$ does not vanish. Hence, we may set $y$ as a new independent variable and consider $(v,\zeta ,J)$ as functions of $y$. Then evolution of $(v,\zeta ,J)$ is described by the following equations  
\begin{eqnarray}
\frac{d v}{d y} &=& 
\nonumber
\frac{2y +\mu k^{1/2}f_2  +\mu^2 k^{1/2} g_2}
{(y^2 -1 - k(v-1)) +\mu k f_1    +\mu^2 k  g_1 },\\
\label{eq_D3}
\frac{d \zeta}{d y} &=& 
\frac{k^{1/2}( \mu f_3   +   \mu^2 g_3 ) }
{(y^2 -1 - k(v-1)) +\mu k f_1    +\mu^2 k  g_1 },\\
\nonumber
\frac{d J}{d y} &=& 
\frac{J_0 k^{1/2}( \mu f_4   +   \mu^2 g_4  )}
{(y^2 -1 - k(v-1)) +\mu k f_1    +\mu^2 k  g_1  }{\rm e}^{-v}.
\end{eqnarray}
We supply these equations by initial conditions corresponding to the point 
$\mathcal{M}^{(2)}$:
\begin{equation}\label{eq_InitCond_3}
v(y^{(2)}) = v^{(2)},\quad \zeta(y^{(2)}) = \zeta^{(2)},\quad J(y^{(2)}) = J^{(2)},
\end{equation}
where $\zeta^{(2)}$, $J^{(2)}$ and $v^{(2)}$ are defined by (\ref{eq_M2}) and (\ref{eq_S2}). 

One represents the solution satisfying (\ref{eq_D3}) and (\ref{eq_InitCond_3}) in the form
\begin{eqnarray}
\nonumber
v(y)&=& v_0(y)  +\mu v_1 (y) +O(\mu^2), \\ 
\label{eq_Expand_3}
\zeta (y)  &=& \zeta_{0} +\mu \zeta_{1}(y)  +\mu^2\zeta_{2} (y) +O(\mu^3),\\ 
\nonumber
J(y) &=& J_{0}+\mu J_{1} (y) +\mu^2 J_{2}(y)  +O(\mu^3),  
\end{eqnarray}
where $v_{0}(y)$ corresponds to the unperturbed orbit $\mathcal{O}(\zeta_{0}, J_{0})$.

Taking into account relations $y=1-kz$ and $v^{(2)}=\frac12 \ln k^{-1}$ and expanding $y^{(2)}$ as
\[ y^{(2)} = y_0^{(2)} + \mu y_1^{(2)} + O(\mu^2),\]
one obtains from (\ref{eq_asymp_z0_D2}) 
\begin{equation}\label{eq_y0_S2}
y_{0}^{(2)}= 1-kz_0(v^{(2)}) = 1 - \frac{1}{2}k^{1/2} + \frac{1}{4}k\ln k^{-1}  + O(k)
\end{equation}
and due to (\ref{eq_z1_S2})
\begin{equation}\label{eq_y1_S2}
y_{1}^{(2)}= - kz_1(v^{(2)}) =  -A(\zeta_{0}) \frac{k^{1/2}}{2J_{0}}
\cdot \Bigr(1 + O(k^{1/2})\Bigr) + O(k^{-5/2}).
\end{equation}

Then, using the Taylor formula with respect to $\mu$, we shift initial conditions (\ref{eq_InitCond_3}) to the point $y=y_0^{(2)}$ as follows:
\begin{align}\label{eq_InitCond_3b}
\nonumber
 & v_0(y_0^{(2)}) = v_{0}^{(2)} ,\quad v_{1}(y_0^{(2)})  =  
- \frac{d v_{0}}{d y}(y_{0}^{(2)})\cdot y_{1}^{(2)},\\
& \zeta_{1}(y_0^{(2)})  = \zeta_{1}^{(2)},\quad
\zeta_{2}(y_0^{(2)}) = \zeta_{2}^{(2)} - \frac{d \zeta_{1}}{d y}(y_{0}^{(2)})\cdot y_{1}^{(2)},\\
\nonumber
& J_{1}(y_0^{(2)})  = J_{1}^{(2)},\quad J_{2}(y_0^{(2)}) = 
J_{2}^{(2)} - \frac{d J_{1}}{d y}(y_{0}^{(2)})\cdot y_{1}^{(2)} .
\end{align} 

Substituting (\ref{eq_Expand_3}) into (\ref{eq_D3}) and expanding with respect to $\mu$, one obtains equations for the $i$-th approximations $\zeta_i$ and $J_i$:
\begin{eqnarray}
\label{Eq:syst_m_D3}
 \frac{d\zeta_{i} }{d y} &=& \mathcal{R}_{i}^{(\zeta)}(y), \\ 
\nonumber
\frac{d J_{i}}{d y } &=& \mathcal{R}_{i}^{(J)}(y).
\end{eqnarray} 

The solution of this system is
\begin{eqnarray}
\label{eq_sol_m_D3}
\zeta_{i}(y) &=&  \zeta_{i}(y_0^{(2)})+\int_{y_{0}^{(2)}}^{y}\mathcal{R}_{i}^{(\zeta)}(s)\, {\rm d}s,\\
\nonumber
J_{i}(v) &=&  J_{i}(y_0^{(2)})+\int_{y_{0}^{(2)}}^{y}\mathcal{R}_{i}^{(J)}(s)\, {\rm d}s.
\end{eqnarray}

For $i=1$ one has
\begin{eqnarray}
\label{eq_rhs_1_D3}
\mathcal{R}_{1}^{(\zeta)}(y) &=& 
k^{1/2}\frac{f_{3}}{y^{2}-1-k (v_{0}(y)-1)},\\
\nonumber
\mathcal{R}_{1}^{(J)}(y) &=& 
J_{0}k^{1/2}\frac{f_{4}}{y^{2}-1-k (v_{0}(y)-1)}{\rm e}^{-v_{0}(y)},
\end{eqnarray}
and for $i=2$
\begin{align}
\label{eq_rhs_2_D3}
&\mathcal{R}_{2}^{(\zeta)}(y) = 
k^{1/2}\frac{f'_{3,\eta}\cdot v_{1}(y) + f'_{3,\zeta}\cdot \zeta_{1}(y)+g_{3}}
{y^{2}-1-k (v_{0}(y)-1)} + 
k^{3/2}\frac{f_{3}\cdot (v_{1}(y) - f_{1})}{(y^{2}-1-k (v_{0}(y)-1))^{2}},\\
\nonumber
&\mathcal{R}_{2}^{(J)}(y) = 
J_{0}k^{1/2}\left[
\frac{(f'_{4,\eta} - f_{4}) \cdot v_{1}(y) + f'_{4,\zeta}\cdot \zeta_{1}(y)+g_{4}}
{y^{2}-1-k (v_{0}(y)-1)} + 
k\frac{f_{4}\cdot (v_{1}(y) - f_{1})}{(y^{2}-1-k (v_{0}(y)-1))^{2}}\right]{\rm e}^{-v_{0}(y)},
\end{align}
where the functions $f_{i}, g_{i}$ are evaluated at the point 
$(\xi, \eta, \zeta) = (k^{-1/2}y, k^{-1} + v_{0}(y) -1, \zeta_{0})$.

Our task is to derive asymptotic formulae for $\zeta_{1,2} (0)$ and $J_{1,2}(0)$. We begin with auxilary asymptotics for $v_0(y)$ and $v_1(y)$.

The function $v_{0}(y)$ satisfies the following equation
\begin{equation}\label{eq_dv0}
\frac{d v_{0}}{d y} = \frac{2y}{y^{2} - 1 - k(v_{0}-1)}.
\end{equation}
As $v_{0}(y)$ corresponds to the unperturbed orbit $\mathcal{O}(\zeta_{0}, J_{0})$  then  due to (\ref{eq_unpertOrbit}) it is described implicitly by an equality
\begin{equation}\label{eq_ev0}
e^{v_0}=k^{-1}(1-y^2+kv_0).
\end{equation}
Hence, in the domain $\mathcal{D}_{3}$ it admits an esimate
\begin{equation}\label{eq_v0_D3}
v_0 (y) =  \ln k^{-1} +\ln (1-y^2) + O(k^{1/2} \ln k^{-1}).
\end{equation}

Applying (\ref{eq_D3}), we obtain the following equation for $v_1(y)$: 
\[
\frac{d v_{1}}{d y} = k\frac{2 y}{(y^2 - 1 - k(v_{0}(y)-1))^2} v_{1} 
+ \mathcal{R}_{1}^{(v)}(y),
\]
where
\[
\mathcal{R}_{1}^{(v)}(y) =
k^{1/2}\frac{ f_{2}}{y^{2}-1-k (v_{0}(y)-1)} - 
k\frac{2 y  f_{1}}{(y^{2}-1-k (v_{0}(y)-1))^{2}} 
\]
and the functions $f_{i}$ are evaluated at the point 
$(\xi, \eta, \zeta) = (k^{-1/2}y, k^{-1} + v_{0}(y) -1, \zeta_{0})$.

The solution of this equation is
\begin{equation} \label{eq_sol_v1}
  v_1 (y) = e^{k\int_{y_{0}^{(2)}}^{y}\frac{2 s}{(s^2 - 1 - k(v_{0}(s)-1))^2}ds}
\left( v_{1}(y_0^{(2)}) + \int_{y_{0}^{(2)}}^{y}  
e^{-k\int_{y_{0}^{(2)}}^{s}\frac{2 p}{(p^2 - 1 - k(v_{0}(p)-1))^2}dp}\cdot
\mathcal{R}_{1}^{(v)}(s)\, ds \right).
\end{equation}

Note that in the domain $\mathcal{D}_{3}$ one has 
\[ \frac{1}{y^2-1} = O(k^{-1/2}), \qquad v_0(y)=O(\ln k^{-1}),\]
\begin{equation}\label{eq_asym_forv1}
  y^2-1-k(v_0-1)=(y^2-1)(1+O(k^{1/2}\ln k^{-1})),
\end{equation}
\[ e^{\pm k\int_{y_{0}^{(2)}}^{y}\frac{2 s}{(s^2 - 1 - k(v_{0}(s)-1))^2}ds} = 1+O(k^{1/2}) , \]
\[ \mathcal{R}_{1}^{(v)}(y)=O(k^{-3/2}) ,
\quad \int_{y_0^{(2)}}^y \mathcal{R}_{1}^{(v)}(s) ds=O(k^{-3/2}) .\]

Taking into account (\ref{eq_InitCond_3b}) and using (\ref{eq_dv0}), (\ref{eq_asym_forv1}), (\ref{eq_y0_S2}), (\ref{eq_y1_S2}),
one concludes from (\ref{eq_sol_v1}) that
\begin{equation}
  \label{eq_v1_D3}
  v_1(y) = - \frac{A(\zeta_0)}{J_0} \Bigl(1+ O
(k^{1/2} \ln k^{-1}) 
\Bigr)+ O(k^{-3/2}).
\end{equation}

Then (\ref{eq_rhs_1_D3}) with (\ref{eq_asym_forv1}), (\ref{eq_ev0}) yield
\begin{eqnarray}
\nonumber
\mathcal{R}_{1}^{(\zeta)}(y) &=& - \frac{k^{-1/2} (\beta_1 +\beta_2 k \zeta_0)}{1-y^2}
(1+O(k^{1/2}\ln k^{-1})) ,\\
\label{eq_as_R_1_D3}
\mathcal{R}_{1}^{(J)}(y) &=& O \left( \frac{J_0k^{-1/2} }{(1-y^2)^2} 
\right) .
\end{eqnarray}
Consequently,
\[ \int_{y_0^{(2)}}^0 \mathcal{R}_{1}^{(\zeta)}(s)\, ds =
k^{-1/2} (\beta_1 + \beta_2 k \zeta_0) 
\left( \frac14 \ln k^{-1} + \frac12 \ln 4 \right)
 +O\left( \ln k^{-1} \right),\] 
\[ \int_{y_0^{(2)}}^0 \mathcal{R}_{1}^{(J)}(s)\, ds = O\left( k^{-1} J_0 \right) .\]
Taking into account  (\ref{eq_InitCond_3b}) and (\ref{eq_as_1_S2}) together with (\ref{eq_sol_m_D3}), one gets
\begin{eqnarray}
\nonumber
\zeta_{1}^{(3)} &=& k^{-3/2} \left[ \frac13 \beta_1 + \beta_2 k\zeta_{0} + \frac{1}{2}\left(\beta_{1}+\beta_{2}k\zeta_{0}\right)k (\ln k^{-1} + \ln 4) - \beta_{1}k + 
O(k^{3/2} \ln k^{-1})\right],\\
\label{eq_as_1_S3_final}
J_{1}^{(3)} &=& A(\zeta_0) 
\Bigl(1+O(k)\Bigr) + O(J_0 k^{-1}) .
\end{eqnarray}

For the second order terms, due to (\ref{eq_sol_m_D3}) and (\ref{eq_InitCond_3b}), we have:
\begin{eqnarray}
\label{eq_sol_m2_D3}
\zeta_{2}(0) &=&  \zeta_{2}^{(2)} - \frac{d \zeta_{1}}{d y}(y_{0}^{(2)})\cdot y_{1}^{(2)}
+\int_{y_{0}^{(2)}}^{0}\mathcal{R}_{2}^{(\zeta)}(s)\, {\rm d}s,\\
\nonumber
J_{2}(0) &=&  J_{2}^{(2)} - \frac{d J_{1}}{d y}(y_{0}^{(2)})\cdot y_{1}^{(2)}
+\int_{y_{0}^{(2)}}^{0}\mathcal{R}_{2}^{(J)}(s)\, {\rm d}s,
\end{eqnarray}
where $\zeta_{2}^{(2)}$ and $J_2^{(2)}$ are defined by (\ref{eq_as_2_S2}),
\[  - \frac{d \zeta_{1}}{d y}(y_{0}^{(2)})\cdot y_{1}^{(2)}
= - \mathcal{R}_{1}^{(\zeta)}(y _{0}^{(2)})\cdot y_{1}^{(2)},
\qquad  - \frac{d J_{1}}{d y}(y_{0}^{(2)})\cdot y_{1}^{(2)}
= - \mathcal{R}_{1}^{(J)}(y_{0}^{(2)})\cdot y_{1}^{(2)} .\]
Formulae (\ref{eq_as_R_1_D3}), (\ref{eq_y0_S2}), (\ref{eq_y1_S2}) imply
\begin{eqnarray}
\label{eq_as_bc_D3}
- \frac{d \zeta_{1}}{d y}(y_{0}^{(2)})\cdot y_{1}^{(2)} &=& 
- A(\zeta_0) \frac{k^{-1/2}}{2J_0}
(\beta_1 +\beta_2 k\zeta_0)
\left( 1+ O (k^{1/2} \ln k^{-1} \right)  ,\\
\nonumber
- \frac{d J_{1}}{d y}(y_{0}^{(2)})\cdot y_{1}^{(2)}&=& O(k^{-1})  .
\end{eqnarray}

Application of (\ref{eq_rhs_2_D3}), (\ref{eq_v0_D3}), (\ref{eq_v1_D3}) leads to the following estimates
\begin{eqnarray}
\label{eq_as_R_2_D3}
\mathcal{R}_{2}^{(\zeta)}(y) &=& 
O\left( \frac{ k^{1/2} }{J_0} \cdot \frac{1}{(1-y^2)^2}
\right) \quad \Rightarrow \int_{y_0^{(2)}}^0 \mathcal{R}_{2}^{(\zeta)}(s)\, ds = O\left( \frac{ 1 }{J_0} \right)  ,\\
\nonumber
\mathcal{R}_{2}^{(J)}(y) &=& O\left( \frac{k^{1/2}}{(1-y^2)^3}  \right) \quad \Rightarrow \int_{y_0^{(2)}}^0 \mathcal{R}_{2}^{(J)}(s)\, ds =  O(k^{-1/2}).
\end{eqnarray}
Therefore, we have the following asymptotics for  the second order approximations
\begin{eqnarray}
\label{eq_as_2_S3_final}
\zeta_{2}^{(3)} &=&  -A(\zeta_0)\frac{k^{-1/2}}{2J_0}\left[  \beta_1 + \beta_2 k\zeta_0 \right] 
\Bigl(1+O(k^{1/2} \ln k^{-1})\Bigr) + O \left( \frac{1}{J_0} \right),\\
\nonumber
J_{2}^{(3)} &=& -A(\zeta_{0})k^{-3/2} \left( \frac{\alpha_{1} }{3} +\alpha_2 k \zeta_0 - \gamma \right)\Bigl(1+O(k^{1/2})\Bigr) + O(k^{-1}).
\end{eqnarray}

Thus, the Poincar\'e map $F(\mu)$ can be represented as
\begin{eqnarray}
\label{eq_F(3)}
F(\mu ) : \left(
\begin{array}{cc}
\zeta_0 \\
J_0 \end{array}
\right) 
\mapsto 
\left(
\begin{array}{cc}
\zeta^{(3)} \\
J^{(3)}\end{array}
\right) =
\left(
\begin{array}{cc}
\zeta_{0} + \mu \zeta_{1}^{(3)}+\mu^2 \zeta_{2}^{(3)} +O(\mu^3) \\
J_{0} +\mu J_{1}^{(3)} +\mu^2 J_{2}^{(3)} +O(\mu^3)\end{array}
\right)
,
\end{eqnarray}
where $\zeta_i^{(3)}$ and $J_i^{(3)}$ are defined by (\ref{eq_as_1_S3_final}) for $i=1$ and (\ref{eq_as_2_S3_final}) for $i=2$.

\section{Fixed points and period-doubling bifurcations}

In this section  we derive conditions on the parameters of the normal form which lead to existence of a fixed point of the Poincar\'e map and its period-doubling bifurcation.
If $(\zeta_0,J_0)\in S_{-}$ is a fixed point then according to (\ref{eq_fixed_point_3}) the following condition should be satisfied:
\[ \zeta_1^{(3)}=0,\quad J_1^{(3)}=0.\]
Taking into account (\ref{eq_as_1_S3_final}), (\ref{eq_def_A}) and definition of the parameter $k$, one may rewrite these conditions as
\begin{align}
\nonumber
&
\frac{\beta_{2}\zeta_{0}}{\ln J_{0}^{-1}}
\left[1 + \frac{\ln\ln J_{0}^{-1}+\ln 4}{2\ln J_{0}^{-1}} \right] +
\beta_{1}\left[\frac{1}{3} + \frac{\ln \ln J_{0}^{-1} + \ln 4 -2}{2\ln J_{0}^{-1}}\right] =
O\left(\frac{\ln\ln J_{0}^{-1}}{\ln^{3/2} J_{0}^{-1}}\right),\\
\label{eq_FxPt}
& - \frac32 \gamma + \alpha_2\zeta_0 -\frac12 \alpha_1 - \gamma_0  
= O( J_{0} \ln^2 J_{0}^{-1} ).
\end{align}

Note that the coefficient $\gamma_0$ in (\ref{eq_FxPt}) is the only one which depends on the ratio, $\delta/\varepsilon$, of parameters of the initial problem (\ref{Eq:initialODE_0}), namely:
\[ \gamma_0 = \varkappa \sigma + \nu, 
\qquad \varkappa = \sqrt{\frac{2}{-D}}F''_{x\delta},
\quad \nu = \sqrt{\frac{2}{-D}}\frac{D'_{x}}{F''_{x^{2}}},
\quad \sigma=\frac{\delta}{\varepsilon}, \]
where the functions $F''_{x\delta}, F''_{x^{2}}, D, D'_{x}$ are evaluated at 
$(x,y,z,\delta)=(0,0,0,0)$.

We solve the first equation with respect to $\zeta_{0}$ and substitute the solution into the second one. Then the second condition (\ref{eq_FxPt}) can be considered as an equation defining $J_{0}$ in terms of the parameter $\sigma$:
\begin{align}
\label{eq_as_FxPt}
&\zeta_{0} = -\frac{\beta_{1}}{3\beta_{2}}\ln J_{0}^{-1}\left(1 + 
\frac{\ln\ln J_{0}^{-1} + \ln 4 - 3}{\ln J_{0}^{-1}}\right) + 
O\left(\frac{\ln \ln J_{0}^{-1}}{\ln^{1/2} J_{0}^{-1}}\right),\\
\nonumber
&\sigma = - \frac{\alpha_{2}\beta_{1}}{3\beta_{2}\varkappa}\ln J_{0}^{-1}\left(1 + \frac{\ln\ln J_{0}^{-1} + \ln 4 - 3}{\ln J_{0}^{-1}}\right) - 
\frac{3\gamma+\alpha_{1}+2\nu}{2\varkappa} + 
O\left(\frac{\ln \ln J_{0}^{-1}}{\ln^{1/2} J_{0}^{-1}}\right). 
\end{align}
We emphasize here that due to $J_{0}\ll 1$ the ratio $\sigma$ satisfies $\sigma\gg 1$.

One may also note that due to smooth dependence of solutions on initial conditions asymptotics  (\ref{eq_as_2_S3_final}) and (\ref{eq_def_A}) being differentiated with respect to $(\zeta_0, J_0)$ remain valid in $U_{\mu}\left(\mathcal{O}(\zeta_{0}, J_{0})\right)$. This leads to
\[ {\rm Tr} \, D_{(\zeta_0, J_0)} F_2 
= - \frac{\sqrt{\pi} k^{-1/2}}{4 J_0} \alpha_2 (\beta_1 + \beta_2 k \zeta_0)
\left( 1 +O(k^{1/2}\ln k^{-1})  \right) +O(k^{-5/2}). \]

Substituting this into (\ref{eq_period_doubling_2}) and taking into account (\ref{eq_FxPt}), one gets a condition for the period-doubling bifurcation of the periodic trajectory corresponding to initial point $(\zeta_{0}, J_{0})$:
\begin{equation*}
1-\mu^2 \frac{\sqrt{\pi} \alpha_{2}\beta_{1} \ln^{1/2}J_{0}^{-1}}{6 J_{0}}
\cdot \left(1 + O\Bigl(\frac{\ln\ln J_{0}^{-1}}{\ln^{1/2}J_{0}^{-1}}\Bigr)   \right)
+ O(\mu^{2}\ln^{-5/2}J_{0}^{-1})=0,
\end{equation*}
where $J_{0}$ satisfies (\ref{eq_as_FxPt}).

We solve this equation with respect to $J_{0}$ and substitute the solution into (\ref{eq_as_FxPt}). Then, taking into account the relation $\varepsilon=\mu^{2}$, one obtains that the first period-doubling bifurcation occurs at 
\begin{align}
\nonumber
&\zeta_{0}^{*} = -\frac{\beta_{1}}{3\beta_{2}}\left[\ln \varepsilon^{-1} +
\frac{1}{2}\ln \ln \varepsilon^{-1} - 
\ln \left(\sqrt{\pi} \frac{\alpha_{2}\beta_{1}{\rm e}^{3}}{24}\right)\right] + 
O\left(\frac{\ln \ln \varepsilon^{-1}}{\ln^{1/2} \varepsilon^{-1}}\right),\\
\label{eq_FPD}
&J_{0}^{*} = \sqrt{\pi}\frac{\alpha_{2}\beta_{1}}{6}\varepsilon\ln^{1/2}\varepsilon^{-1}
\left(1 + O\left(\frac{\ln \ln \varepsilon^{-1}}{\ln \varepsilon^{-1}}\right)\right),\\
\nonumber
&\delta^{*} = 
	\begin{aligned}[t]
		- \frac{\alpha_{2}\beta_{1}}{3\beta_{2}\varkappa}\varepsilon
		\biggl[\ln \varepsilon^{-1} +
		\frac{1}{2}\ln \ln \varepsilon^{-1} - 
		&\ln \left(\sqrt{\pi} \frac{\alpha_{2}\beta_{1}{\rm e}^{3}}{24}\right)\biggr] 
		- \\
		&-\frac{3\gamma+\alpha_{1}+2\nu}{2\varkappa}\varepsilon + 
		O\left(\varepsilon\frac{\ln \ln \varepsilon^{-1}}{\ln^{1/2} \varepsilon^{-1}}				\right).
	\end{aligned}
\end{align}

It is to be noted that (\ref{eq_as_FxPt}), (\ref{eq_FPD}) hold true provided $\varkappa \ne 0$, $\alpha_2 \ne 0$, $\beta_{1,2} \ne 0$. We apply (\ref{eq_Coef2}), (\ref{eq_Coef3}), (\ref{eq_Coef}) to obtain:
\[
\beta_1 = \frac{\sqrt{2}}{F''_{x^2}\sqrt{-D}}
\left(  G'_{1x} \left| \begin{array}{ccc}
F''_{x^2} & G''_{1x^2} & G''_{2x^2} \\
0 & G'_{1x} & G'_{2x}\\
F'_{y} & G'_{1y} & G'_{2y}
\end{array}
\right| + G'_{2x} \left| \begin{array}{ccc}
F''_{x^2} & G''_{1x^2} & G''_{2x^2} \\
0 & G'_{1x} & G'_{2x}\\
F'_{z} & G'_{1z} & G'_{2z}
\end{array}
\right| 
\right),
\]

\[
\alpha_2 = \frac{\sqrt{2}}{(-D)^{5/2}}
\left(  F'_{y} \left| \begin{array}{ccc}
F''_{x^2} & F''_{xy} & F''_{xz} \\
0 & F_{y} & F_{z}\\
G'_{1x} & G'_{1y} & G'_{1z}
\end{array}
\right| + F'_{z} \left| \begin{array}{ccc}
F''_{x^2} & F''_{xy} & F''_{xz} \\
0 & F'_{y} & F'_{z}\\
G'_{2x} & G'_{2y} & G'_{2z}
\end{array}
\right| 
\right),
\]

\[ \beta_2 = \frac{\sqrt{2}}{(-D)^{3/2}} \left| \begin{array}{ccc}
0 & F'_y & F'_z \\
G'_{1x} & G'_{1y} & G'_{1z} \\
G'_{2x} & G'_{2y} & G'_{2z}
\end{array}
\right| .
\]
If all coefficients $\varkappa$, $\alpha_{2}$, $\beta_{1,2}$ vanish then fixed point does not exist for small values of $J_{0}$. In other cases one needs to perform further asymptotic analysis to obtain conditions for existence of a periodic orbit and its period-doubling bifurcation.

One may also remark that the distance between the fixed point and the fold is  
\[ 
\rho = \delta F''_{x\delta}\sqrt{\frac{F'^{2}_{y} + F'^{2}_{z}}{(F'_{y}F''_{xz})^{2} + F''^{2}_{x^{2}}(F'^{2}_{y} + F'^{2}_{z})}} + O(\delta^{2}).
\]
Thus, the cascade of period-doubling bifurcations occurs when the equilibrium is not very close to the fold, but situated at a distance of the order 
$O(\varepsilon \ln \varepsilon^{-1})$.

\section{Example: the FitzHugh-Nagumo system}

We apply our results to the FitzHugh-Nagumo system (\ref{eq_FHN}).
Let
\[ \delta = 1-a,\qquad a+x=:x,\quad y+a-\frac{a^3}{3} =:y.\]
Then the system takes the form:
\begin{eqnarray}
\nonumber
\varepsilon \dot x&=&-\frac13 x^3 + x^2 (1-\delta ) +x(2\delta -\delta^2)-y-z \,, \\
\label{eq_FHN_NF}
\dot y&=&  x \,, \\
\nonumber
\dot z &=&  x-z \, 
\end{eqnarray}
and the functions $F$, $G_1$ and $G_2$ are
\[ F(x,y,z,\delta) =-\frac13 x^3 + x^2 (1-\delta ) +x(2\delta -\delta^2)-y-z ,\quad G_1 (x,y,z,\delta) =x
,\quad G_2 (x,y,z,\delta) = x-z.\]
Note that the conditions (\ref{eq_Cond_RHS1}), (\ref{eq_Cond_RHS2}) and (\ref{eq_Cond_RHS}) are satisfied.

Make a rescaling of parameters
\[ \varepsilon = \mu^2 , \quad  \delta = \mu^2 \sigma, \]
and introduce new variables $(\xi, \eta, \zeta)$ by
\begin{equation*}
x = \mu \xi \,, \quad
y = \frac{\mu^2}{2} (\eta-\zeta)\, , \quad
z = \frac{\mu^2}{2} (\eta +\zeta) \, .
\end{equation*}

Then the inverse change gives 
$\mu=  \sqrt{\varepsilon }$, $\sigma= \delta  \varepsilon^{-1}$ and
\begin{equation*}
\xi = \frac{1}{\mu }\, x\,, \quad
\eta = \frac{1}{\mu^2}( y+z) \,, \quad
\zeta = \frac{1}{\mu^2}\left( - y+z \right)\, .
\end{equation*}

In terms of these variable the FitzHugh-Nagumo system takes the form
\begin{eqnarray}
\nonumber
\xi'&=&\xi^2-\eta 
+ \mu \left( 2\sigma \xi - \frac{1}{3}\xi^3 \right) - \mu^2 \sigma\xi^2 -\mu^3\sigma^2\xi \,, \\
\label{eq_FHN_NF}
\eta'&=&  2\xi + \mu \left( -\frac12 \eta - \frac{1}{2} \zeta \right) \,, \\
\nonumber
\zeta' &=&   \mu \left( -\frac{1}{2}\eta - \frac12 \zeta\right) \, .
\end{eqnarray}

Then 
\[ \gamma_0=2\sigma, \quad \gamma=-\frac13 ,\quad
\alpha_1=\alpha_2=\beta_1=\beta_2=-\frac12,\quad
\varkappa = 2,\quad \nu = 0 \]
and condition (\ref{eq_FPD}) reads:
\[ \delta = 
\frac{1}{12}\varepsilon
\left[\ln \varepsilon^{-1} +
\frac{1}{2}\ln \ln \varepsilon^{-1} - 
\ln \left(\frac{\sqrt{\pi}{\rm e}^{3}}{96}\right)\right] + 
\frac{3}{8}\varepsilon + 
O\left(\varepsilon\frac{\ln \ln \varepsilon^{-1}}{\ln^{1/2} \varepsilon^{-1}}\right).
\]

We compare our results with numerical data obtained by M. Zaks and found sufficiently good agreement.

\begin{table}[h]
\begin{tabular}{llll}
$\varepsilon$& $a_{num}$		& $a_{asym}$		& $a_{num}-a_{asym}$\\
1.e-2		&0.99092058501692	&0.99094938062714	&2.879561021731e-5\\
1.e-4  	&0.99986818929447	&0.99986822927480	&3.99803325e-8\\
1.e-6  	&0.99999828100195	&0.99999828163419	&6.322363e-10\\
1.e-8  	&0.99999997885167	&0.99999997885883	&7.1557e-12\\
1.e-10		&0.99999999974920	&0.99999999974928	&7.28e-14\\
1.e-12		&0.99999999999710	&0.99999999999710	&7.e-16
\end{tabular}
\caption{Comparison of numerical and asymptotic results: values of the parameter $a$ at the first period-doubling for several values of $\varepsilon$.}
\end{table}

\end{document}